\documentclass[graybox]{svmult}

\usepackage{mathptmx}       
\usepackage{helvet}         
\usepackage{courier}        
\usepackage{type1cm}        
\usepackage{makeidx}         
\usepackage{graphicx}        
\usepackage{multicol}        
\usepackage{multirow}
\usepackage[bottom]{footmisc}
\usepackage{bm}
\usepackage{marginnote}

\makeindex             

\usepackage{benstyle_book}

\begin{document}

\title*{Compressed sensing approaches for polynomial approximation of high-dimensional functions}
\author{Ben Adcock, Simone Brugiapaglia and Clayton G.\ Webster}
\institute{Ben Adcock \at Simon Fraser University, Burnaby, BC, Canada 
\email{ben\_adcock@sfu.ca}
\and
Simone Brugiapaglia \at Simon Fraser University, Burnaby, BC, Canada \email{simone\_brugiapaglia@sfu.ca}
\and 
Clayton G.\ Webster \at University of Tennessee and Oak Ridge National Lab, Oak Ridge, TN, USA \email{webstercg@math.utk.edu, webstercg@ornl.gov}
}
%
%
\maketitle

\abstract{In recent years, the use of sparse recovery techniques in the approximation of high-dimensional functions has garnered increasing interest.  In this work we present a survey of recent progress in this emerging topic.  Our main focus is on the computation of polynomial approximations of high-dimensional functions on $d$-dimensional hypercubes.  We show that smooth, multivariate functions possess expansions in orthogonal polynomial bases that are not only approximately sparse, but possess a particular type of structured sparsity defined by so-called lower sets.  This structure can be exploited via the use of weighted $\ell^1$ minimization techniques, and, as we demonstrate, doing so leads to sample complexity estimates that are at most logarithmically dependent on the dimension $d$.  Hence the curse of dimensionality -- the bane of high-dimensional approximation -- is mitigated to a significant extent.  We also discuss several practical issues, including unknown noise (due to truncation or numerical error), and highlight a number of open problems and challenges. 
}

\section{Introduction}\label{s:introduction}
The approximation of high-dimensional functions is a fundamental difficulty
in a large number of fields, including 
neutron, tomographic and magnetic resonance image reconstruction, uncertainty quantification (UQ), optimal control and parameter identification for engineering and science applications.  In addition, this problem 
naturally arises in computational solutions to kinetic plasma physics equations, the many-body Schr\"odinger equation, Dirac and Maxwell equations for molecular electronic structures and nuclear dynamic computations, options pricing equations in mathematical finance, Fokker-Planck and fluid dynamics equations for complex fluids, turbulent flow, quantum dynamics, molecular life sciences, and nonlocal mechanics.
The subject of intensive research over the last half-century, high-dimensional approximation is made challenging by the \textit{curse of dimensionality}, a phrase coined by Bellman \cite{bellmann}.  Loosely speaking, this refers to the tendency of 
na\"ive approaches to exhibit exponential blow-up in complexity with increasing dimension.  Progress is possible, however, by placing restrictions on the class of functions to be approximated; for example, smoothness, anisotropy, sparsity, and compressibility.  Well-known algorithms such as sparse grids \cite{Webster:2007tq,bunggrieb,Nobile:2008wf,Nobile:2008uc}, which are specifically designed to capture this behaviour, can mitigate the curse of dimensionality to a substantial extent. 

While successful, however, such approaches typically require strong \textit{a priori} knowledge of the functions being approximated, e.g.\ the parameters of the anisotropic behaviour, or costly adaptive implementations to estimate the anisotropy during the approximation process.  The efficient approximation of high-dimensional functions in the absence of such knowledge remains a significant challenge. 

In this chapter, we consider new methods for high-dimensional approximation based on the techniques of compressed sensing.  Compressed sensing is an appealing approach for reconstructing signals from underdetermined systems, i.e.\ with far smaller number of measurements compared to the signal 
length \cite{CRT06,Donoho06}. This approach has emerged in the last half a dozen years as an alternative to more classical approximation schemes for high-dimensional functions, with the aim being to overcome some of the limitations mentioned above.  Under natural the sparsity or compressibility assumptions, it enjoys a significant improvement in sample complexity 
over traditional methods such as discrete least-squares, projection, and interpolation 
\cite{Gunzburger:2014hi,Gunzburger:2016gt}. 
Our intention in this chapter is to both present an overview of existing work in this area, focusing particularly on the mitigation of the curse of dimensionality, and to highlight existing open problems and challenges.

\subsection{Compressed sensing for high-dimensional approximation}
Compressed sensing asserts that a vector $\bm{x} \in \bbC^{n}$ possessing a $k$-sparse representation in a fixed orthonormal basis can be recovered from a number of suitably-chosen measurements $m$ that is linear in $k$ and logarithmic in the ambient dimension $n$.  In practice, recovery can be achieved via a number of different approaches, including convex optimization ($\ell^1$ minimization), greedy or thresholding algorithms.

Let $f : D \rightarrow \bbC$ be a function, where $D \subseteq \bbR^d$ is a domain in $d \gg 1$ dimensions.  In order to apply compressed sensing techniques to approximate $f$, we must first address the following three questions:
\begin{enumerate}
\item[(i)] In which orthonormal system of functions $\{ \phi_{i} \}^{n}_{i=1}$ does $f$ have an approximately sparse representation?
\item[(ii)] Given suitable assumptions on $f$ (e.g.\ smoothness) how fast does the best $k$-term approximation error decay?
\item[(iii)] Given such a system $\{ \phi_{i} \}^{n}_{i=1}$, what are suitable measurements to take of $f$?
\end{enumerate}
The concern of this chapter is the approximation of smooth functions, and as such we will use orthonormal bases consisting of multivariate orthogonal polynomials.  In answer to (i) and (ii) in \S \ref{s:sparsepoly} we discuss why this choice leads to approximate sparse representations for functions with suitable smoothness, and characterize the best $k$-term approximation error in terms of certain regularity conditions.  As we note in \S \ref{ss:regularity}, practical examples of such functions include parameter maps of many different types of parametric PDEs.

For sampling, we evaluate $f$ at a set of points $\bm{z}_1,\ldots,\bm{z}_m \in D$.  This approach is simple and particularly well-suited in practical problems.  In UQ, for example, it is commonly referred to as a \textit{nonintrusive} approach \cite{SMUQ} or \textit{stochastic collocation} \cite{NarayanZhouCCP}.  More complicated measurement procedures -- for instance, \textit{intrusive} procedures such as inner products with respect to a set of basis functions -- are often impractical or even infeasible, since, for example, they require computation of high-dimensional integrals.  The results presented in \S \ref{s:CShighdim} identify appropriate (random) choices of the sample points $\{ \bm{z}_i \}^{m}_{i=1}$ and bounds for the number of measurements $m$ under which $f$ can be stably and robustly recovered from the data $\{ f(\bm{z}_i) \}^{m}_{i=1}$.

\subsection{Structured sparsity}
The approximation of high-dimensional functions using polynomials differs from standard compressed sensing in several key ways.  Standard compressed sensing exploits sparsity of the finite vector of coefficients $\bm{c} \in \bbC^n$ of a finite-dimensional signal $\bm{x} \in \bbC^n$.  However, polynomial coefficients of smooth functions typically possess more detailed structure than just sparsity.  Loosely speaking, coefficients corresponding to low polynomial orders tend to be larger than coefficients corresponding to higher orders.  This raises several questions:
\begin{enumerate}
\item[(iv)] What is a reasonable \textit{structured sparsity} model for polynomial coefficients of high-dimensional functions?
\item[(v)] How can such structured sparsity be exploited in the reconstruction procedure, and by how much does doing this reduce the number of measurements required?
\end{enumerate}
In \S \ref{ss:lowersets} it is shown that high-dimensional functions can be approximated with quasi-optimal rates of convergence by $k$-term polynomial expansions with coefficients lying in so-called \textit{lower} sets of multi-indices.  As we discuss, sparsity in lower sets is a type of structured sparsity, and in \S \ref{s:CShighdim} we show how it can be exploited by replacing the classical $\ell^1$ regularizer by a suitable weighted $\ell^1$-norm.  Growing weights penalize high degree polynomial coefficients, and when chosen appropriately, they act to promote lower set structure.  In \S \ref{ss:nonuniform} nonuniform recovery techniques are used to identify a suitable choice of weights.  This choice of weights is then adopted in \S \ref{ss:uniform} to establish quasi-optimal uniform recovery guarantees for compressed sensing of polynomial expansions using weighted $\ell^1$ minimization.  

The effect of this weighted procedure is a substantially improved recovery guarantee over the case of unweighted $\ell^1$ minimization; specifically, a measurement condition that is only logarithmically dependent on the dimension $d$ and polynomially-dependent on the sparsity $k$.  Hence the curse of dimensionality is almost completely avoided.  
As we note in \S \ref{ss:comparison}, these polynomial rates of growth in $k$ agree with the best known recovery guarantees for oracle least-squares estimators.


\subsection{Dealing with infinity}\label{ss:infinity}

Another way in which the approximation of high-dimensional functions differs from  standard compressed sensing is that functions typically have infinite (as opposed to finite) expansions in orthogonal polynomial bases.  In order to apply compressed sensing techniques, this expansion must be truncated in a suitable way.  This leads to the following questions:
\begin{enumerate}
\item[(vi)] What is a suitable truncation of the infinite expansion?
\item[(vii)] How does the corresponding truncation error affect the overall reconstruction?
\end{enumerate}  
In \S \ref{s:CShighdim} a truncation strategy -- corresponding to a hyperbolic cross index set -- is proposed based on the lower set structure.  The issue of truncation error (question (vii)) presents some technical issues, both theoretical and practical, since this error is usually unknown \textit{a priori}.  In \S \ref{ss:functions} we discuss a means to overcome these issues via a slightly modified optimization problem.  Besides doing so, another benefit of the approach developed therein is that it yields approximations to $f$ that also interpolate at the sample points $\{ \bm{z}_i \}^{m}_{i=1}$; a desirable property for certain applications.  Furthermore, the results given in \S \ref{ss:functions} also address the robustness of the recovery to unknown errors in the measurements.  This is a quite common phenomenon in applications, since function samples are often the result of (inexact) numerical computations.

\subsection{Main results}\label{ss:main_result}
We now summarize our main results.  In order to keep the presentation brief, in this chapter we limit ourselves to functions defined on the unit hypercube $D = (-1,1)^d$ and consider expansions in orthonormal polynomial bases $\{ \phi_{\bm{i}} \}_{\bm{i} \in \bbN^d_0}$ of Chebyshev or Legendre type.  We note in passing, however, that many of our results apply immediately (or extend straightforwardly) to more general systems of functions.  See \S \ref{s:conclusions} for some further discussion.

Let $\nu$ be the probability measure under which the basis $\{ \phi_{\bm{i}} \}_{\bm{i} \in \bbN^d_0}$ is orthonormal.  Our main result is as follows:

\thm{
\label{t:intro_thm}
Let $k \in \bbN$, $0 < \epsilon < 1$, $\{ \phi_{\bm{i}} \}_{\bm{i} \in \bbN^d_0}$ be the orthonormal Chebyshev or Legendre basis on $D = (-1,1)^d$, $\Lambda = \Lambda^{\mathrm{HC}}_{k}$ be the hyperbolic cross of index $k$ and define weights $\bm{u} = (u_{\bm{i}})_{\bm{i} \in \bbN^d_0}$, where $u_{\bm{i}} = \nm{\phi_{\bm{i}}}_{L^\infty}$.  Suppose that
\bes{
m \gtrsim k^{\gamma}  \left ( \log^2(k) \min \left \{ d + \log(k) , \log(2d) \log(k) \right \} + \log(k) \log(\log(k)/\epsilon) \right ),
}
where $\gamma = \frac{\log(3)}{\log(2)}$ or $\gamma = 2$ for Chebyshev or Legendre polynomials respectively, and draw $\bm{z}_{1},\ldots,\bm{z}_{m} \in D$ independently according to $\nu$.  Then with probability at least $1-\epsilon$ the following holds.  For any $f \in L^2_{\nu}(D) \cap L^{\infty}(D)$ satisfying
\be{
\label{f_tail}
\nm{ f - \sum_{\bm{i} \in \Lambda} c_{\bm{i}} \phi_{\bm{i}} }_{L^\infty} \leq \eta,
}
for some known $\eta \geq 0$, it is possible to compute, via solving a $\ell^1_{\bm{u}}$ minimization problem of size $m \times n$ where $n = | \Lambda |$, an approximation $\tilde{f}$ from the samples $\bm{y} = (f(\bm{z}_j))^{m}_{j=1}$ that satisfies
\be{
\label{intro_err}
\| f - \tilde{f} \|_{L^2_{\nu}} \lesssim \frac{\sigma_{k,L}(\bm{c})_{1,\bm{u}} }{k^{\gamma/2}} + \eta,\qquad \nm{f - \tilde{f}}_{L^\infty} \lesssim \sigma_{k,L}(\bm{c})_{1,\bm{u}} + k^{\gamma/2} \eta .
}
Here $\bm{c}$ are the coefficients of $f$ in the basis $\{ \phi_{\bm{i}} \}_{\bm{i} \in \bbN^d_0}$ and $\sigma_{s,L}(\bm{c})_{1,\bm{u}}$ is the $\ell^1_{\bm{u}}$-norm error of the best approximation of $\bm{c}$ by a vector that is $k$-sparse and lower.
}

%

Note that the condition \R{f_tail} is strong, since it assumes an \textit{a priori} upper bound on the expansion error is available.  Such a condition is unlikely to be met in practice.  In \S \ref{ss:functions} we discuss recovery results for general $f$ without such \textit{a priori} bounds.

\subsection{Existing literature}\label{ss:literature}

The first results on compressed sensing with orthogonal polynomials in the one-dimensional setting appeared in \cite{Rauhut}, based on earlier work in  sparse trigonometric expansions \cite{RauhutTrigPoly}.  This was extended to the higher-dimensional setting in \cite{YanGuoXui_l1UQ}.  Weighted $\ell^1$ minimization was studied in \cite{RauhutWardWeighted}, and recovery guarantees given in terms of so-called weighted sparsity.  However, this does not lead straightforwardly to explicit measurement conditions for quasi-best $k$-term approximation.  The works \cite{AdcockCSFunInterp,ChkifaDownwardsCS} introduced new guarantees for weighted $\ell^1$ minimization of nonuniform and uniform types respectively, leading to optimal sample complexity estimates for recovering high-dimensional functions using $k$-term approximations in lower sets.  Theorem \ref{t:intro_thm} is based on results in \cite{ChkifaDownwardsCS}.  Relevant approaches to compressed sensing in infinite dimensions have also been considered in \cite{AdcockCSFunInterp,Adcockl1Pointwise,BAACHGSCS,BrugiapagliaThesis,BNMP15,TraonmilinGribonvalRIP}

Applications of compressed sensing to UQ, specifically the computation of polynomial chaos expansions of parametric PDEs, can be found in \cite{BouchotEtAlMultilevel,DoostanOwhadiSparse,MathelinGallivanCSPDErandom,PengHamptonDoostantweighted,RS14,KarniadakisUQCS} and references therein. 
Throughout this chapter we use random sampling from the orthogonality measure of the polynomial basis.  We do this for its simplicity and the theoretical optimality of the recovery guarantees in terms of the dimension $d$.  Other strategies, which typically seek a smaller error or lower polynomial factor of $k$ in the sample complexity, have been considered in \cite{GuoEtAlRandomizedQuad,HamptonDoostanCSPCE,JakemanEtAlChristoffel,NarayanJakemanZhouChristoffelLS,NarayanZhouCCP,TangIaccarino,XuZhouSparseDeterministic}.  Working towards a similar end, various approaches have also been considered to learn a better sparsity basis \cite{JakemanEtAl_l1Enhance,YangEtAlEnhancingRotations} or to use additional gradient samples \cite{PengHampDoos15}.  In this chapter, we focus on fixed bases of Chebyshev or Legendre polynomials in the unit cube.  For results in $\bbR^d$ using Hermite polynomials, see \cite{GuoEtAlRandomizedQuad,HamptonDoostanCSPCE,NarayanJakemanZhouChristoffelLS}.

In some scenarios, a suitable lower set may be known in advance or be computed via an adaptive search.  In this case, least-squares methods may be suitable.  A series of works have studied the sample complexity of such approaches in the context of high-dimensional polynomial approximation \cite{ChkifaEtAl,DavenportEtAlLeastSquares,MiglioratiCohenOptimal,CMN15,HamptonDoostantCohMot,MiglioratiThesis,MiglioratiJAT,MiglioratiNobileLowDisc,MiglioratiEtAlFoCM,NarayanJakemanZhouChristoffelLS}.  We review a number of these results in \S \ref{ss:comparison}.

\section{Sparse polynomial approximation of high-dimensional functions}\label{s:sparsepoly}


\subsection{Setup and notation}\label{ss:setup}
We first require some notation.  For the remainder of this chapter, $D = (-1,1)^d$ will be the $d$-dimensional unit cube.  The vector $\bm{z} = (z_1,\ldots,z_d)$ will denote the variable in $D$ and $\bm{i} = (i_1,\ldots,i_d) \in \bbN^d_0$ will be a multi-index.  Let $\nu^{(1)},\ldots,\nu^{(d)}$ be probability measures on the unit interval $(-1,1)$.  We consider the tensor product probability measure $\nu$ on $D$ given by $\nu = \nu^{(1)} \otimes \cdots \otimes \nu^{(d)}$.  Let $\{ \phi^{(k)}_{i} \}^{\infty}_{i=0}$ be an orthonormal polynomial basis of $L^2_{\nu^{(k)}}(-1,1)$ and define the corresponding tensor product orthonormal basis $\{ \phi_{\bm{i}} \}_{\bm{i} \in \bbN^d_0}$ of $L^2_{\nu}(D)$ by
\bes{
\phi_{\bm{i}} = \phi^{(1)}_{i_1} \otimes \cdots \otimes \phi^{(d)}_{i_d},\qquad \bm{i} = (i_1,\ldots,i_d) \in \bbN^d_0.
}
We let $\nm{\cdot}_{L^2_{\nu}}$ and $\ip{\cdot}{\cdot}_{L^2_{\nu}}$ denote the norm and inner product on $L^2_{\nu}(D)$ respectively.

Let $f \in L^2_{\nu}(D) \cap L^{\infty}(D)$ be the function to be approximated, and write
\be{
\label{f_expansion}
f = \sum_{\bm{i} \in \bbN^d_0} c_{\bm{i}} \phi_{\bm{i}},
}
where $c_{\bm{i}} = \ip{f}{\phi_{\bm{i}}}_{L^2_{\nu}}$ are the coefficients of $f$ in the basis $\{ \phi_{\bm{i}} \}_{\bm{i} \in \bbN^d_0}$.  We define
\bes{
\bm{c} = \left ( c_{\bm{i}} \right )_{\bm{i} \in \bbN^d_0} \in \ell^2(\bbN^d_0),
}
to be the infinite vector of coefficients in this basis.

\examp{
\label{ex:Jacobi}
Our main example will be Chebyshev or Legendre polynomials.  In one dimension, these are orthogonal polynomials with respect to the weight functions
\eas{
& \D \nu = \frac12 \D z \quad \mbox{(Legendre)},\qquad \D \nu = \frac{1}{\pi \sqrt{1-z^2}} \D z\quad \mbox{(Chebyshev)},
}
respectively.  For simplicity, we will consider only tensor products of the same types of polynomials in each coordinate.  The corresponding tensor product measures on $D$ are consequently defined as
\eas{
& \D \nu = 2^{-d} \D \bm{z} \quad \mbox{(Legendre)},\qquad \D \nu = \prod^{d}_{j=1} \frac{1}{\pi \sqrt{1-z^2_j}} \D \bm{z}\quad \mbox{(Chebyshev)}.
}
We note also that many of the results presented below extend to more general families of orthogonal polynomials, e.g., Jacobi polynomials (see Remark \ref{r:Jacobi}).
}

As discussed in \S \ref{ss:infinity}, it is necessary to truncate the infinite expansion \R{f_expansion} to a finite one.  Throughout, we let $\Lambda \subset \bbN^d_0$ be a subset of size $| \Lambda | = n$, and define the truncated expansion
\bes{
f_{\Lambda} = \sum_{\bm{i} \in \Lambda} c_{\bm{i}} \phi_{\bm{i}}.
}
We write $\bm{c}_{\Lambda}$ for the finite vector of coefficients with multi-indices in $\Lambda$.  Whenever necessary, we will assume an ordering $\bm{i}_1,\ldots,\bm{i}_n$ of the multi-indices in $\Lambda$, so that
\bes{
f_{\Lambda} = \sum^{n}_{j=1} c_{\bm{i}_j} \phi_{\bm{i}_j},\qquad \bm{c}_{\Lambda} = ( c_{\bm{i}_j} )^{n}_{j=1} \in \bbC^n.
}
We will adopt the usual convention and view $\bm{c}_{\Lambda}$ interchangeably as a vector in $\bbC^n$ and as an element of $\ell^2(\bbN^d_0)$ whose entries corresponding to indices $\bm{i} \notin \Lambda$ are zero.

\subsection{Regularity and best $k$-term approximation}\label{ss:regularity}

In the high-dimensional setting, we assume the regularity of $f$ is such that the complex continuation of $f$, represented as the map $f:\mathbb{C}^d\to \mathbb{C}$, is a  holomorphic function on $\mathbb{C}^d$.  In addition,  for $1 \leq k \leq n$, we let
\bes{
\Sigma_{k} = \left \{ \bm{c} \in \ell^2(\bbN^d_0): | \mathrm{supp}(\bm{c}) | \leq k \right \},
}
be the set of $k$-sparse vectors, and
\bes{
\sigma_{k}(\bm{c})_1 = \inf_{\bm{d} \in \Sigma_k} \| \bm{c}-\bm{d} \|_{1},
}
be the error of the best $k$-term approximation of $\bm{c}$, measured in the $\ell^1$ norm.

Recently, for smooth functions as described above, sparse recovery of the polynomial expansion \eqref{f_expansion} with the use of compressed sensing has shown tremendous promise.  
However, this approach requires a {small uniform bound} of the underlying basis, given by
\begin{align*}
\Theta = \sup_{{\bm i}\in \Lambda} \|{ \phi}_{\bm i}\|_{L^{\infty}({D})},
\end{align*}
as the sample complexity $m$ required to recover the best $k$-term approximation (up to a multiplicative constant) scales with the following bound (see, e.g., \cite{FoucartRauhutCSbook})
\begin{align}
\label{intro_cond}
m\gtrsim \Theta^2 k \times \text{log factors}.
\end{align}
This poses a challenge for many multivariate polynomial approximation strategies  
as $\Theta$ is prohibitively large in high dimensions. In particular, for $d$-dimensional problems, 
$\Theta = 2^{d/2}$ for Chebyshev systems and so-called preconditioned 
Legendre systems \cite{Rauhut}. Moreover, when using the standard Legendre expansion, the 
theoretical number of samples can exceed the cardinality of the polynomial subspace, unless 
the subspace a priori excludes all terms of high total order (see, e.g., \cite{HamptonDoostanCSPCE,YanGuoXui_l1UQ}). 
Therefore, the advantages of sparse polynomial recovery methods,  
coming from reduced sample complexity,  are eventually overcome by the curse of dimensionality, 
in that such techniques require at least as many samples as traditional sparse interpolation techniques
in high dimensions \cite{Gunzburger:2014hi,Nobile:2008uc,Nobile:2008wf}.  Nevertheless, in the next section we 
describe a common characteristic of the polynomial space spanned by the best $k$-terms, that we will exploit to  overcome the curse of dimensionality in the sample complexity bound \eqref{intro_cond}.
As such, our work also provides a fair comparison with existing numerical polynomial approaches 
in high dimensions \cite{BNTT14,CCDS13,ChkifaEtAl,CCS13,Tran:2015C}.

\subsection{Lower sets and structured sparsity}\label{ss:lowersets}

{
In many engineering and science applications, the target functions, despite being high-dimensional, are smooth and often characterized by a rapidly decaying polynomial expansion, whose most important coefficients are {of} {low order} {\cite{ChkifaEtAlBreaking,CD15,CDS11,HoangSchwab12}.} 
In such situations, the quest for finding the approximation containing the largest $k$ terms
can be restricted to polynomial spaces associated with {\em lower (also known as downward closed or monotone)} sets.  These are defined as follows:
%
\defn{
A set $S \subseteq \bbN^{d}_{0}$ is lower if, whenever $\bm{i} = (i_1,\ldots,i_d) \in S$ and $\bm{i'} = (i'_1,\ldots,i'_d) \in \bbN^d_0$ satisfies $i'_j \leq i_j$ for all $j=1,\ldots,d$, then $\bm{i'} \in S$.
}

{The practicality of downward closed sets is mainly computational, and has been 
demonstrated in different approaches such as quasi-optimal strategies, Taylor expansion, interpolation methods, 
and discrete least-squares (see
\cite{BNTT14,CCDS13,ChkifaEtAl,CCS13, 
ChkifaEtAlBreaking,ChkifaDownwardsCS,CohenDeVoreSchwabFoCM,
CohenDeVoreSchwabRegularity,MiglioratiJAT,MiglioratiEtAlFoCM,Stoyanov:2016ci,Tran:2015C} and references therein).  
For instance, in the context of parametric PDEs, it was shown in \cite{ChkifaEtAlBreaking} that for a large class of smooth differential operators, with a certain type of anisotropic dependence on ${\bm z}$, the solution map ${\bm z}\mapsto f({\bm z})$ can be approximated by its best $k$-term expansions associated with {index sets of cardinality $k$}, resulting in algebraic rates $k^{-\alpha},\, \alpha>0$ in the uniform and/or mean average sense. The same rates are preserved with index sets that are lower.  
In addition, such lower sets of cardinality $k$ also enable the equivalence property 
$\|\cdot\|_{L^2_{\nu}(D)}\leq \|\cdot\|_{L^\infty} \leq k^{\gamma}\|\cdot\|_{L^2_{\nu}(D)}$ in arbitrary dimensions $d$ with, e.g.,  
$\gamma=2$ for the uniform measure and $\gamma=\frac{\log3}{\log2}$ for Chebyshev measure.
%
%

Rather than best $k$-term approximation, we now consider best $k$-term approximation in a lower set.  Hence, we replace $\Sigma_{k}$ with
\bes{
\Sigma_{k,L} = \left \{ \bm{c} \in \ell^2(\bbN^d_0) : | \mathrm{supp}(\bm{c}) | \leq k,\ \mbox{$\mathrm{supp}(\bm{c})$ is lower} \right \},
}
and $\sigma_{k}(\bm{c})_1$ with the quantity
\be{
\label{sigma_s_lower}
\sigma_{k,L}(\bm{c})_{1,\bm{w}} = \inf_{\bm{d} \in \Sigma_{k,L}} \| \bm{c}-\bm{d} \|_{1,\bm{w}}.
}
Here $\bm{w} = ( w_{\bm{i}} )_{\bm{i} \in \bbN^d_0}$ is a sequence of positive weights and $\nm{\bm{c}}_{1,\bm{w}} = \sum_{\bm{i} \in \bbN^d_0} w_{\bm{i}} | c_{\bm{i}} |$ is the norm on $\ell^1_{\bm{w}}(\bbN^d_0)$.  

\rem{
\label{r:structuredsparsity}
Sparsity in lower sets is an example of a so-called \textit{structured sparsity} model.  Specifically, $\Sigma_{k,L}$ is the subset of $\Sigma_{k}$ corresponding to the union of all $k$-dimensional subspaces defined by lower sets:
\bes{
\Sigma_{k,L} \equiv \bigcup_{\substack{|S| = k \\ \mbox{\small $S$ lower}}} \left \{ \bm{c}: \supp(\bm{c}) \subseteq S \right \} \subset  \bigcup_{|S| = k} \left \{ \bm{c}: \supp(\bm{c}) \subseteq S \right \} \equiv \Sigma_{k} .
}
Structured sparsity models have been studied extensively in compressed sensing (see, e.g.,  \cite{BaranuikModelCS,BlumensathUnionSubspace,EldarDuarteCSReview,DuarteEldarStructuredCS,TraonmilinGribonvalRIP} and references therein).  There are a variety of general approaches for exploiting such structure, including greedy and iterative methods (see, for example, \cite{BaranuikModelCS}) and convex relaxations \cite{TraonmilinGribonvalRIP}.  A difficulty with lower set sparsity is that projections onto $\Sigma_{k,L}$ cannot be easily computed \cite{ChkifaDownwardsCS}, unlike the case of $\Sigma_{k}$.  Therefore, in this chapter we shall opt for a different approach based on $\ell^1_{\bm{w}}$ minimization with suitably-chosen weights $\bm{w}$.  See \S \ref{s:conclusions} for some further discussion.
}

\section{Compressed sensing for multivariate polynomial approximation}\label{s:CShighdim}

Having introduced tensor orthogonal polynomials as a good basis for obtaining (structured) sparse representation of smooth, multivariate functions, we now turn our attention to computing quasi-optimal approximations of such a function $f$ from the measurements $\{ f(\bm{z}_i) \}^{m}_{i=1}$.

It is first necessary to choose the sampling points $\bm{z}_1,\ldots,\bm{z}_m$.  From now on, following an approach that has become standard in compressed sensing \cite{FoucartRauhutCSbook}, we shall assume that these points are drawn randomly and independently according to the probability measure $\nu$.  We remark in passing that this may not be the best choice in practice.  However, such an approach yields recovery guarantees with measurement conditions that are essentially independent of $d$, thus mitigating the curse of dimensionality.  In \S \ref{s:conclusions} we briefly discuss other strategies for choosing these points which may convey some practical advantages.


\subsection{Exploiting lower set-structured sparsity}

Let $\bm{c} \in \ell^2(\bbN^d_0)$ be the infinite vector of coefficients of a function $f \in L^2_{\nu}(D)$.  Suppose that $\Lambda \subset \bbN^d_0$, $| \Lambda | =n$ and notice that
\be{
\label{underdet_system}
\bm{y} = A \bm{c}_{\Lambda} + \bm{e}_{\Lambda},
}
where $\bm{y} \in \bbC^m$ and $A \in \bbC^{m \times n}$ are the finite vector and matrix given by
\be{
\label{y_A_def}
\bm{y} = \frac{1}{\sqrt{m}} \left ( f(\bm{z}_j) \right )^{m}_{j=1},\qquad A = \frac{1}{\sqrt{m}} \left ( \phi_{\bm{i}_k}(\bm{z}_j) \right )^{m,n}_{j,k=1},
}
respectively, and
\be{
\label{e_Lambda_def}
\bm{e}_{\Lambda} = \frac{1}{\sqrt{m}} \left ( f(\bm{z}_j) - f_{\Lambda}(\bm{z}_j) \right )^{m}_{j=1} = \frac{1}{\sqrt{m}} \left ( \sum_{\bm{i} \notin \Lambda} c_{\bm{i}} \phi_{\bm{i}}(\bm{z}_j) \right )^{m}_{j=1},
}
is the vector of remainder terms corresponding to the coefficients $c_{\bm{i}}$ with indices outside $\Lambda$.  Our aim is to approximate $\bm{c}$ up to an error depending on $\sigma_{k,L}(\bm{c})_{1,\bm{w}}$, i.e.\ its best $k$-term approximation in a lower set (see \R{sigma_s_lower}).  In order for this to be possible, it is necessary to choose $\Lambda$ so that it contains all lower sets of cardinality $k$.  A straightforward choice is to make $\Lambda$ exactly equal to the union of all such sets, which transpires to be precisely the hyperbolic cross index set with index $k$.  That is,
\be{
\label{HCdef}
\bigcup_{\substack{|S| = k \\ \mbox{\small $S$ lower}}} S = \left \{ \bm{i} = (i_1,\ldots,i_d) \in \bbN^d_0 : \prod^{d}_{j=1} (i_j+1) \leq k \right \} = \Lambda^{\mathrm{HC}}_{k}.
}
It is interesting to note that this union is a finite set, due to the lower set assumption.  Had one not enforced this additional property, the union would be infinite and equal to the whole space $\bbN^d_0$.

We shall assume that $\Lambda = \Lambda^{\mathrm{HC}}_{k}$ from now on.  For later results, it will be useful to know the cardinality of this set.  While an exact formula in terms of $k$ and $d$ is unknown, there are a variety of different upper and lower bounds.  In particular, we shall make use of the following result:
\be{
\label{HCsize}
n = \left | \Lambda^{\mathrm{HC}}_{k} \right  | \leq \min \left \{ 2 k^3 4^d , \E^2 k^{2 + \log_2(d)} \right \}.
}
See \cite[Thm.\ 3.7]{ChernovDungHCCardinality} and \cite[Thm.\ 4.9]{KuhnEtAlApproxMixed} respectively.

With this in hand, we now wish to obtain a solution $\hat{\bm{c}}_{\Lambda}$ of \R{underdet_system} which approximates $\bm{c}_{\Lambda}$, and therefore $\bm{c}$ due to the choice of $\Lambda$, up to an error determined by its best approximation in a lower set of size $k$.  We shall do this by weighted $\ell^1$ minimization.  Let $\bm{w} = \left ( w_{\bm{i}} \right )_{\bm{i} \in \Lambda}$ be a vector of positive weights and consider 
the problem
\be{
\label{weighted_l1_min}
\min_{\bm{d} \in \bbC^n} \nm{\bm{d}}_{1,\bm{w}}\ \mbox{s.t.}\ \nm{\bm{y} - A \bm{d} }_{2} \leq \eta,
}
where $\nm{\bm{d}}_{1,\bm{w}} = \sum^{n}_{j=1} w_{\bm{i}_j} | d_{\bm{i}_j} |$ is the weighted $\ell^1$-norm and $\eta \geq 0$ is a parameter that will be chosen later.  Since the weights $\bm{w}$ are positive we shall without loss of generality assume that
\bes{
w_{\bm{i}} \geq 1,\quad \forall \bm{i}.
}
Our choice of these weights is based on the desire to exploit the lower set structure.  Indeed, since lower sets inherently penalize higher indices, it is reasonable (and will turn out to be the case) that appropriate choices of increasing weights will promote this type of structure.

For simplicity, we shall assume for the moment that $\eta$ is chosen so that
\be{
\label{tail_bound}
\eta \geq \nm{\bm{e}_{\Lambda}}_{2}.
}
In particular, this implies that the exact vector $\bm{c}_{\Lambda}$ is a feasible point of the problem \R{weighted_l1_min}.  As was already mentioned in \S \ref{ss:main_result}, this assumption is a strong one, and is unreasonable for practical scenarios where good \textit{a priori} estimates on the expansion error $f - f_{\Lambda}$ are hard to obtain.  In \S \ref{ss:functions} we address the removal of this condition.

\subsection{Choosing the optimization weights: nonuniform recovery}\label{ss:nonuniform}

Our first task is to determine a good choice of optimization weights.  For this,  techniques from nonuniform recovery\footnote{By nonuniform recovery, we mean results that guarantee recovery of a fixed vector $\bm{c}_{\Lambda}$ from a single realization of the random matrix $A$.  Conversely, uniform recovery results consider recovery of all sparse (or structured sparse) vectors from a single realization of $A$.  See, for example, \cite{FoucartRauhutCSbook} for further discussion.} are particularly useful.  

At this stage it is convenient to define the following.  First, for a vector of weights $\bm{w}$ and a subset $S$ we let
\be{
\label{weighted_cardinality}
| S |_{\bm{w}} = \sum_{\bm{i} \in S} w^2_{\bm{i}},
}
be the \textit{weighted} cardinality of $S$.  Second, for the orthonormal basis $\{ \phi_{\bm{i}} \}_{\bm{i} \in \bbN^d_0}$ we define the \textit{intrinsic} weights $\bm{u} = \left ( u_{\bm{i}} \right )_{\bm{i}}$ as
\be{
\label{intrinsic_weights}
u_{\bm{i}} = \nm{\phi_{\bm{i}}}_{L^\infty}.
}
Note that $u_{\bm{i}} = \nm{\phi_{\bm{i}}}_{L^\infty} \geq\nm{\phi_{\bm{i}}}_{L^2_{\nu}} =1$ since $\nu$ is a probability measure.  With this in hand, we now have the following result (see \cite[Thm.\ 6.1]{AdcockCSFunInterp}):

\thm{
\label{t:nonuniformCS}
Let $\Lambda \subset \bbN^d_0$ with $| \Lambda | = n $, $0 < \epsilon < \E^{-1}$, $\eta \geq 0$, $\bm{w} = (w_{\bm{i}} )_{\bm{i} \in \Lambda}$ be a set of weights, $\bm{c} \in \ell^2(\bbN^d_0)$ and $S \subseteq \Lambda$, $S \neq \emptyset$, be any fixed set.  Draw $\bm{z}_{1},\ldots,\bm{z}_m$ independently according to the measure $\nu$, let $A$, $\bm{y}$ and $\bm{e}_{\Lambda}$ be as in \R{y_A_def} and \R{e_Lambda_def} respectively and suppose that $\eta$ satisfies \R{tail_bound}.  Then, with probability at least $1-\epsilon$, any minimizer $\hat{\bm{c}}_{\Lambda}$ of \R{weighted_l1_min} satisfies
\be{
\label{nonunif_error_bound}
\nm{\bm{c} - \hat{\bm{c}}_{\Lambda} }_2 \lesssim \lambda \sqrt{| S |_{\bm{w} }} \left ( \eta + \| \bm{c} - \bm{c}_{\Lambda} \|_{1,\bm{u}} \right ) + \nm{\bm{c} - \bm{c}_{S} }_{1,\bm{w}},
}
provided
\be{
\label{nonunif_meas_cond}
m \gtrsim \left ( | S |_{\bm{u}} + \max_{\bm{i} \in \Lambda \backslash S} \{ u^2_{\bm{i}} / w^2_{\bm{i}} \}  | S |_{\bm{w}} \right )  L,
}
where $\lambda = 1 + \frac{\sqrt{\log(\epsilon^{-1})}}{\log(2n \sqrt{| S |_{\bm{w}} } ) }$ and $L = \log(\epsilon^{-1})  \log \left(2 n \sqrt{| S |_{\bm{w}} } \right )$.
}

Suppose for simplicity that $\bm{c}$ were exactly sparse and let $S = \mathrm{supp}(\bm{c})$ and $\eta = 0$.  Then this result asserts exact recovery of $\bm{c}$, provided the measurement condition \R{nonunif_meas_cond} holds.  Ignoring the log factor $L$, this condition is determined by
\be{
\label{Mdef}
\cM (S; \bm{u} , \bm{w}) =  | S |_{\bm{u}} + \max_{\bm{i} \in \Lambda \backslash S} \{ u^2_{\bm{i}} / w^2_{\bm{i}} \} | S |_{\bm{w}} .
}
The first term is the weighted cardinality of $S$ with respect to the intrinsic weights $\bm{u}$, and is independent of the choice of optimization weights $\bm{w}$.  The second term depends on these weights, but the possibly large size of $| S |_{\bm{w}}$ is compensated by the factor $\max_{\bm{i} \in \Lambda \backslash S} \{ u^2_{\bm{i}} / w^2_{\bm{i}} \}$.

Seeking to minimize $\cM (S; \bm{u} , \bm{w})$, it is natural to choose the weights $\bm{w}$ so that the second term in \R{Mdef} is equal to the first.  This is easily achieved by the choice
\be{
\label{weights_choice}
w_{\bm{i}} = u_{\bm{i}},\quad \forall \bm{i},
}
with the resulting measurement condition being simply
\be{
\label{sample_complexity}
m \gtrsim | S |_{\bm{u}}  \log(\epsilon^{-1})  \log(2n \sqrt{|S|_{\bm{u}}}).
}
From now on, we primarily consider the weights \R{weights_choice}.

\rem{
Theorem \ref{t:nonuniformCS} is a nonuniform recovery guarantee for weighted $\ell^1$ minimization.  Its proof uses the well-known golfing scheme \cite{GrossGolfing}, following similar arguments to those given in \cite{BAACHGSCS,Candes_Plan} for unweighted $\ell^1$ minimization.  Unlike the results in \cite{BAACHGSCS,Candes_Plan}, however, it gives a measurement condition in terms of a fixed set $S$, rather than the sparsity $k$ (or weighted sparsity).  In other words, no sparsity (or structured sparsity) model is required at this stage.  Such an approach was first pursued in \cite{BigotBlockCS} in the context of block sampling in compressed sensing.  See also \cite{AdcockChunParallel}.
}

\subsection{Comparison with oracle estimators}\label{ss:comparison}
As noted above, the condition \R{sample_complexity} does not require $S$ to be a lower set.  In \S \ref{ss:oracle} we shall use this property in order to estimate $|S|_{\bm{u}}$ in terms of the sparsity $k$.  First, however, it is informative to compare \R{sample_complexity} to the measurement condition of an oracle estimator.  Suppose that the set $S$ were known.  Then a standard estimator for $\bm{c}$ is the least-squares solution
\be{
\check{\bm{c}}_{S} = (A_{S})^{\dag} \bm{y},
} 
where $A_{S} \in \bbC^{m \times |S|}$ is the matrix formed from the columns of $A$ with indices belonging to $S$ and $\dag$ denotes the pseudoinverse.  Stable and robust recovery via this estimator follows if the matrix $A_{S}$ is well-conditioned.  For this, one has the following well-known result:

\prop{
Let $0 < \delta, \epsilon < 1$, $S \subset \bbN^d_0$, $|S| = k$ and suppose that $m$ satisfies
\be{
\label{LS_meas_cond}
m \gtrsim \delta^{-2}  | S |_{\bm{u}}  \log(2k\epsilon^{-1}).
}
Draw $\bm{z}_1,\ldots,\bm{z}_m$ independently according to the measure $\nu$ and let $A$ be as in \R{y_A_def}.  Then, with probability at least $1-\epsilon$, the matrix $A_{S}$ satisfies
\bes{
\nm{ (A_S)^* A_S - I }_{2} \leq \delta,
}
where $I \in \bbC^{k \times k}$ is the identity matrix and $\nm{\cdot}_2$ is the matrix $2$-norm.
}

See, for example, \cite[Lem.\ 8.2]{AdcockCSFunInterp}.  Besides the log factor, \R{LS_meas_cond} is the same sufficient condition as \R{sample_complexity}.  Thus the weighted $\ell^1$ minimization estimator $\hat{\bm{c}}_{\Lambda}$ with weights $\bm{w} = \bm{u}$ requires roughly the same measurement condition as the oracle least-squares estimator.  Of course, the former requires no \textit{a priori} knowledge of $S$.

\rem{
In fact, one may prove a slightly sharper estimate than \R{LS_meas_cond} where $|S|_{\bm{u}}$ is replaced by the quantity
\be{
\label{Christoffel_fn}
\sup_{\bm{z} \in D} \sum_{\bm{i} \in S} |\phi_{\bm{i}}(\bm{z}) |^2 .
}
See, for example, \cite{DavenportEtAlLeastSquares}.  Note that $\sum_{\bm{i} \in S} |\phi_{\bm{i}}(\bm{z})|^2$ is the so-called Christoffel function of the subspace spanned by the functions $\{ \phi_{\bm{i}} \}_{\bm{i} \in S}$.  However, \R{Christoffel_fn} coincides with $|S|_{\bm{u}}$ whenever the polynomials $\phi_{\bm{i}}$ achieve their absolute maxima at the same point in $D$.  This is the case for any Jacobi polynomials whenever the parameters satisfy $\max \{ \alpha,\beta \} \geq -1/2$ \cite[Thm.\ 7.32.1]{SzegoOrthPolys}; in particular, Legendre and Chebyshev polynomials (see Example \ref{ex:Jacobi}), and tensor products thereof.
}

\subsection{Sample complexity for lower sets}\label{ss:oracle}
The measurement condition \R{sample_complexity} determines the sample complexity in terms of the weighted cardinality $|S|_{\bm{u}}$ of the set $S$.  When a structured sparsity model is applied to $S$ -- in particular, lower set sparsity -- one may derive estimates for $|S|_{\bm{u}}$ in terms of just the cardinality $k = |S|$ and the dimension $d$.

\lem{
\label{l:lower_christoffel}
Let $2 \leq k \leq 2^{d+1}$.  If $\{\phi_{\bm{i}} \}_{\bm{i}\in\bbN^d_0}$ is the tensor Chebyshev basis then
\bes{
k^{\log(3)/\log(2)} / 3 \leq \max \left \{ |S|_{\bm{u}} : S \subset \bbN^d_0,\ |S| \leq k, \mbox{$S$ lower} \right \} \leq k^{\log(3)/\log(2)} ,
}
where $|S|_{\bm{u}}$ and $\bm{u}$ are as in \R{weighted_cardinality} and \R{intrinsic_weights} respectively.  If $\{\phi_{\bm{i}} \}_{\bm{i}\in\bbN^d_0}$ is the tensor Legendre basis then
\bes{
k^2 /4 \leq \max \left \{ |S|_{\bm{u}} : S \subset \bbN^d_0,\ |S| \leq k, \mbox{$S$ lower} \right \} \leq k^2.
}
Moreover, the upper estimates hold for all $k \geq 2$.
}

See \cite[Lem.\ 3.7]{ChkifaEtAl}.  With this in hand, we now have the following result:

\thm{
\label{t:nonuniformChebLeg}
Consider the setup in Theorem \ref{t:nonuniformCS} with $k \geq 2$, $\Lambda = \Lambda^{\mathrm{HC}}_{k}$ the hyperbolic cross \R{HCdef}, weights $\bm{w} = \bm{u}$ and $\{ \phi_{\bm{i}} \}_{\bm{i} \in \bbN^d_0}$ the tensor Legendre or Chebyshev basis.  Then any minimizer $\hat{\bm{c}}_{\Lambda}$ of \R{weighted_l1_min} with weights $\bm{w} = \bm{u}$ satisfies
\bes{
\nm{\bm{c}-\hat{\bm{c}}_{\Lambda} }_{2} \lesssim \lambda k^{\gamma/2} \left ( \eta + \nm{\bm{c}-\bm{c}_{\Lambda}}_{1,\bm{u}} \right ) + \sigma_{k,L}(\bm{c})_{1,\bm{u}},
}
with probability at least $1-\epsilon$, provided
\bes{
m \gtrsim k^{\gamma} \log(\epsilon^{-1}) \min \left \{ d + \log(k) , \log(2d) \log(k) \right \},
}
where $\lambda = 1 + \frac{\sqrt{\log(\epsilon^{-1})}}{\log(k)}$ and where $\gamma = \log(3)/\log(2)$ or $\gamma = 2$ in the Chebyshev or Legendre case respectively.
}
\prf{
Let $S \subset \bbN^d_0$, $|S| \leq k$ be a lower set such that $\nm{\bm{c} - \bm{c}_{S}}_{1,\bm{u}} = \sigma_{k,L}(\bm{c})_{1,\bm{u}}$.  By Lemma \ref{l:lower_christoffel} we have $| S |_{\bm{u}} \leq k^{\gamma}$.  We now apply Theorem \ref{t:nonuniformCS} with $\bm{w} = \bm{u}$, and use this result and the bound \R{HCsize} for $n = | \Lambda^{\mathrm{HC}}_{k}|$.
\qed
}

\rem{
It is worth noting that the lower set assumption drastically reduces the sample complexity.  Indeed, for the case of Chebyshev polynomials one has
\bes{
\max \left \{ |S|_{\bm{u}} : S \subset \bbN^d_0,\ |S| \leq k\right \} = 2^d k.
}
In other words, in the absence of the lower set condition, one can potentially suffer exponential blow-up with dimension $d$.  Note that this result follows straightforwardly from the explicit expression for the weights $\bm{u}$ in this case: namely,
\be{
\label{Chebu}
u_{\bm{i}} = 2^{\nm{\bm{i}}_0/2},
}
where $\nm{\bm{i}}_0 = | \left \{ j : i_{j} \neq 0 \right \} |$ for $\bm{i} = (i_1,\ldots,i_d) \in \bbN^d_0$ (see, for example, \cite{AdcockCSFunInterp}).  On the other hand, for Legendre polynomials the corresponding quantity is infinite, since in this case the weights
\be{
\label{Legu}
u_{\bm{i}} =  \prod^{d}_{j=1} \sqrt{2 i_j+1},
}
are unbounded.  Moreover, even if $S$ is constrained to lie in the hyperbolic cross $\Lambda = \Lambda^{\mathrm{HC}}_{k}$, one still has a worst-case estimate that is polynomially-large in $k$ \cite{ChkifaDownwardsCS}.
}

\rem{
\label{r:Jacobi}
We have considered only tensor Legendre and Chebyshev polynomial bases.  However, Theorem \ref{t:nonuniformChebLeg} readily extends to other types of orthogonal polynomials.  All that is required is an upper bound for
\bes{
\max \left \{ |S|_{\bm{u}} : S \subset \bbN^d_0,\ |S| \leq k, \mbox{$S$ lower} \right \},
}
in terms of the sparsity $k$.  For example, suppose that $\{ \phi_{\bm{i}} \}_{\bm{i} \in \bbN^d_0}$ is the tensor ultraspherical polynomial basis, corresponding to the measure
\bes{
\D \nu =  (c_{\alpha})^d\prod^{d}_{j=1} (1-z^2_j)^{\alpha} \D \bm{z},\qquad c_{\alpha} = \left ( \int^{1}_{-1} (1-z^2)^{\alpha} \D z \right )^{-1}.
}
If the parameter $\alpha$ satisfies $2 \alpha + 1 \in \bbN$ then \cite[Thm.\ 8]{MiglioratiJAT} gives that 
\bes{
\max \left \{ |S|_{\bm{u}} : S \subset \bbN^d_0,\ |S| \leq k, \mbox{$S$ lower} \right \} \leq k^{2 \alpha + 2}.
}
This result includes the Legendre case ($\alpha = 0$) given in Lemma \ref{l:lower_christoffel}, as well as the case of Chebyshev polynomials of the second kind ($\alpha = 1/2$).  A similar result also holds for tensor Jacobi polynomials for parameters $\alpha,\beta \in \bbN_0$ (see \cite[Thm.\ 9]{MiglioratiJAT}).
}

\subsection{Quasi-optimal approximation: uniform recovery}\label{ss:uniform}
As is typical of a nonuniform recovery guarantee, the error bound in Theorem \ref{t:nonuniformChebLeg} has the limitation that it relates the $\ell^2$-norm of the error with the best $k$-term, lower approximation error in the $\ell^1_{\bm{u}}$-norm.  To obtain stronger estimates we now consider uniform recovery techniques.

We first require an extension of the standard Restricted Isometry Property (RIP) to the case of sparsity in lower sets.  To this end, for $k \in \bbN$ we now define the quantity
\be{
\label{sk_def}
s(k) = \max \left \{ |S|_{\bm{u}} : S \subset \bbN^d_0,\ |S| \leq k, \mbox{$S$ lower} \right \}.
}
The following extension of the RIP was introduced in \cite{ChkifaDownwardsCS}:

\defn{
A matrix $A \in \bbC^{m \times n}$ has the lower Restricted Isometry Property (lower RIP) of order $k$ if there exists as $0 < \delta < 1$ such that
\bes{
\left ( 1-\delta \right ) \nm{\bm{c}}^2_{2} \leq \nm{A \bm{c}}^2_2 \leq \left ( 1 + \delta \right) \nm{\bm{c}}^2_2,\quad \forall \bm{c} \in \bbC^n,\ |\supp(\bm{c}) |_{\bm{u}} \leq s(k).
}
If $\delta = \delta_{k,L}$ is the smallest constant such that this holds, then $\delta_{k,L}$ is the $k^{\rth}$ lower Restricted Isometry Constant (lower RIC) of $A$.
}

We shall use the lower RIP to establish stable and robust recovery.  For this, we first note that the lower RIP implies a suitable version of the robust Null Space Property (see \cite[Prop.\ 4.4]{ChkifaDownwardsCS}):

\lem{
\label{LRIPimpliesNSP}
Let $k \geq 2$ and $A \in \bbC^{m \times n}$ satisfy the lower RIP of order $\alpha k$ with constant
\bes{
\delta = \delta_{\alpha k, L} < 1/5,
}
where $\alpha = 2$ if the weights $\bm{u}$ arise from the tensor Legendre basis and $\alpha = 3$ if the weights arise from the tensor Chebyshev basis.  Then for any $S \subseteq \Lambda^{\mathrm{HC}}_{s}$ with $| S |_{\bm{u}} \leq s(k)$ and any $\bm{d} \in \bbC^n$,
\bes{
\nm{\bm{d}_{S}}_{2} \leq \frac{\rho}{\sqrt{s(k)}} \nm{\bm{d}_{S^c}}_{1,\bm{u}} + \tau \nm{A \bm{d}}_{2},
}
where $\rho = \frac{4 \delta}{1-\delta} < 1$ and $\tau = \frac{\sqrt{1+\delta}}{1-\delta}$.
}

With this in hand, we now establish conditions under which the lower RIP holds for matrices $A$ defined in \R{y_A_def}.  The following result was shown in \cite{ChkifaDownwardsCS}:

\thm{
\label{t:lowerRIPm}
Fix $0 < \epsilon < 1$, $0 < \delta < 1/13$, let $\{ \phi_{\bm{i}} \}_{\bm{i} \in \bbN^d_0}$ be as in \S \ref{ss:setup} and $\bm{u}$ be as in \R{intrinsic_weights} and suppose that
\bes{
m \gtrsim \frac{s(k)}{\delta^2} L,
}
where $s(k)$ is as in \R{sk_def} and
\bes{
L = \log \left ( \frac{s(k)}{\delta^2} \right ) \left ( \frac{1}{\delta^4} \log \left ( 2 \frac{s(k)}{\delta^2} \log \left ( \frac{s(k)}{\delta^2} \right ) \right ) \log(2n) + \frac{1}{\delta} \log \left ( \frac{1}{\gamma \delta} \log \left ( \frac{k(s)}{\delta^2} \right ) \right ) \right ).
}
Draw $\bm{z}_{1},\ldots,\bm{z}_{m}$ independently according to $\nu$ and let $A \in \bbC^{m \times n}$ be as in \R{y_A_def}.  Then with probability at least $1-\epsilon$ the matrix $A$ satisfies the lower RIP of order $k$ with constant $\delta_{k,L} \leq 13 \delta$.
}

Combining this with the previous lemma now gives the following uniform recovery guarantee:

\thm{
\label{t:unifquasiopt}
Let $0 < \epsilon < 1$, $k \geq 2$ and 
\be{
\label{uniformsamplecomplexity}
m \asymp k^{\gamma}  L,
}
where $\gamma = \log(3)/\log(2)$ or $\gamma = 2$ in the tensor Chebyshev or tensor Legendre cases respectively and
\be{
\label{uniformlogfactor}
L = \left ( \log^2(k) \min \left \{ d + \log(k) , \log(2d) \log(k) \right \} + \log(k) \log(\log(k)/\epsilon) \right ).
}
Let $\Lambda = \Lambda^{\mathrm{HC}}_{k}$ be the hyperbolic cross index set, $\{ \phi_{\bm{i}} \}_{\bm{i} \in \bbN^d_0}$ be the tensor Legendre or Chebyshev polynomial basis and draw $\bm{z}_{1},\ldots,\bm{z}_{m}$ independently according to the corresponding measure $\nu$.  Then with probability at least $1-\epsilon$ the following holds.  For any $f \in L^2(D) \cap L^\infty(D)$ the approximation
\bes{
\tilde{f} = \sum_{\bm{i} \in \Lambda} \hat{c}_{\bm{i}} \phi_{\bm{i}},
}
where $\bm{\hat{c}}_{\Lambda} = (\hat{c}_{\bm{i}})_{\bm{i} \in \Lambda}$ is a solution of \R{weighted_l1_min} with $A$, $\bm{y}$ and $\eta$ given by \R{y_A_def} and \R{tail_bound} respectively and weights $\bm{w} = \bm{u}$, satisfies
\be{
\label{ferr_unif}
\nm{ f - \tilde{f} }_{L^\infty} \leq \nm{ \bm{c} - \hat{\bm{c}}_{\Lambda} }_{1,\bm{u}} \lesssim \sigma_{k,L}(\bm{c})_{1,\bm{u}} + k^{\gamma/2} \eta,
}
and
\be{
\label{ferr_L2}
\nm{ f - \tilde{f} }_{L^2_{\nu}} = \nm{ \bm{c} - \hat{\bm{c}}_{\Lambda} }_{2}  \lesssim \frac{\sigma_{k,L}(\bm{c})_{1,\bm{u}}}{k^{\gamma/2}} + \eta,
}
where $\bm{c} \in \ell^2(\bbN^d_0)$ are the coefficients of $f$ in the basis  $\{ \phi_{\bm{i}} \}_{\bm{i} \in \bbN^d_0}$.
}
\prf{
Let $\alpha = 2$ or $\alpha = 3$ in the Legendre or Chebyshev case respectively.  Condition \R{uniformsamplecomplexity}, Lemma \ref{l:lower_christoffel} and Theorem \ref{t:lowerRIPm} imply that the matrix $A$ satisfies the lower RIP of order $\alpha k$ with constant $\delta_{\alpha k , L } \leq 1/6 < 1/5$.  Now let $S$ be a lower set of cardinality $|S| = k$ such that
\be{
\label{Sdef1}
\nm{\bm{c} - \bm{c}_{S}}_{1,\bm{u}} = \sigma_{k,L}(\bm{c})_{1,\bm{u}},
}
set $\bm{d} = \bm{c}_{\Lambda} - \hat{\bm{c}}_{\Lambda}$ and $T = \Lambda \backslash S$.  Note that
\eas{
\nm{\bm{d}_{T}}_{1,\bm{u}} &\leq \nm{\bm{c}_{T}}_{1,\bm{u}} + \nm{\hat{\bm{c}}_{T}}_{1,\bm{u}} 
\\
& = 2 \nm{\bm{c}_{T}}_{1,\bm{u}} + \nm{ \bm{c}_{S}}_{1,\bm{u}} +  \nm{\hat{\bm{c}}_{T}}_{1,\bm{u}}  - \nm{\bm{c}_{\Lambda}}_{1,\bm{u}}
\\
& \leq 2 \nm{\bm{c}_{T}}_{1,\bm{u}}  +  \nm{ \bm{d}_{S}}_{1,\bm{u}} + \nm{\hat{\bm{c}}_{\Lambda}}_{1,\bm{u}}  - \nm{\bm{c}_{\Lambda}}_{1,\bm{u}}
 \leq 2 \sigma_{k,L}(\bm{c})_{1,\bm{u}} +  \nm{ \bm{d}_{S}}_{1,\bm{u}},
}
since $\hat{\bm{c}}_{\Lambda}$ is a solution of \R{weighted_l1_min} and $\bm{c}_{\Lambda}$ is feasible for \R{weighted_l1_min} due to the choice of $\eta$. By Lemma \ref{LRIPimpliesNSP} we have
\eas{
\nm{\bm{d}_{T}}_{1,\bm{u}} &\leq 2 \sigma_{k,L}(\bm{c})_{1,\bm{u}} + \sqrt{s(k)} \nm{ \bm{d}_{S}}_2
 \leq 2 \sigma_{k,L}(\bm{c})_{1,\bm{u}} + \rho \nm{\bm{d}_T}_{1,\bm{u}} + \tau \sqrt{s(k)} \nm{A \bm{d}}_{2},
}
where $\rho \leq 4/5$ and $\tau \leq \sqrt{42}/5$.  Therefore
\bes{
\nm{\bm{d}_{T}}_{1,\bm{u}} \lesssim \sigma_{k,L}(\bm{c})_{1,\bm{u}} + \sqrt{s(k)} \nm{A \bm{d}}_{2} \lesssim  \sigma_{k,L}(\bm{c})_{1,\bm{u}} + \sqrt{s(k)} \eta,
}
where in the second step we use the fact that $\bm{d} = \bm{c}_{\Lambda} - \hat{\bm{c}}_{\Lambda}$ is the difference of two vectors that are both feasible for \R{weighted_l1_min}.  Using this bound and Lemma \ref{LRIPimpliesNSP}  again gives
\be{
\label{z1u_bound}
\nm{\bm{d}}_{1,\bm{u}} \lesssim \sigma_{k,L}(\bm{c})_{1,\bm{u}} + \sqrt{s(k)} \eta,
}
and since $\bm{c} - \hat{\bm{c}}_{\Lambda} = \bm{d} + \bm{c} - \bm{c}_{\Lambda}$, we deduce that
\be{
\label{cchat_bound}
\nm{\bm{c} - \hat{\bm{c}}_{\Lambda}}_{1,\bm{u}} \leq \nm{\bm{d}}_{1,\bm{u}} + \nm{\bm{c} - \bm{c}_{\Lambda}}_{1,\bm{u}} \lesssim \sigma_{k,L}(\bm{c})_{1,\bm{u}} + \sqrt{s(k)} \eta.
}
Due to the definition of the weights $\bm{u}$ we have 
$
\nm{f - \tilde{f}}_{L^\infty} \leq \nm{\bm{c} - \hat{\bm{c}}_{\Lambda}}_{1,\bm{u}},
$
and therefore, after noting that $s(k) \lesssim k^{\gamma}$ (see Lemma \ref{l:lower_christoffel}) we obtain the first estimate \R{ferr_unif}.  For the second estimate let $S$ be such that
\bes{
\nm{\bm{c} - \bm{c}_{S}}_{2} = \min \left \{ \nm{\bm{c} - \bm{d}}_{2} : | \supp(\bm{d}) |_{\bm{u}} \leq s(k) \right \},
}
and set $T = S^c$.  Let $\bm{d} = \bm{c}- \hat{\bm{c}}_{\Lambda}$ and write
$
\nm{\bm{d}}_{2} \leq \nm{\bm{d}_{S}}_{2} + \nm{\bm{d}_{T}}_{2}.
$
Via a weighted Stechkin estimate \cite[Thm.\ 3.2]{RauhutWardWeighted} we have
$
\nm{\bm{d}_{T}}_{2} \leq \frac{1}{\sqrt{ s(k) - \nm{\bm{u}}_{\infty}}} \nm{\bm{d}}_{1,\bm{u}}.
$
For tensor Chebyshev and Legendre polynomials, one has $\nm{\bm{u}}_{\infty} \leq \frac34 s(k)$ (see \cite[Lem.\ 4.1]{ChkifaDownwardsCS}), and therefore
$
\nm{\bm{d}_{T}}_{2} \lesssim \frac{1}{\sqrt{s(k)}} \nm{\bm{d}}_{1,\bm{u}}.
$
We now apply Lemma \ref{LRIPimpliesNSP} to deduce that
$
\nm{\bm{d}}_{2} \lesssim \frac{1}{\sqrt{s(k)}} \nm{\bm{d}}_{1,\bm{u}} + \eta.
$
Recall that $s(k) \gtrsim k^\gamma$ due to Lemma \ref{l:lower_christoffel}.  Hence \R{cchat_bound} now gives
$
\nm{\bm{d}}_{2} \lesssim \frac{\sigma_{k,L}(\bm{c})_{1,\bm{u}} }{k^{\gamma/2}} + \eta,
$
as required.
\qed
}

For the Legendre and Chebyshev cases, Theorem \ref{t:unifquasiopt} proves recovery with quasi-optimal $k$-term rates of approximation subject to the same measurement condition (up to log factors) as the oracle least-squares estimator.  In particular, the sample complexity is polynomial in $k$ and at most logarithmic in the dimension $d$, thus mitigating the curse of dimensionality to a substantial extent.  We remark in passing that this result can be extended to general Jacobi polynomials (recall Remark \ref{r:Jacobi}).  Furthermore, the dependence on $d$ can be removed altogether by considering special classes of lower sets, known as \textit{anchored} sets \cite{CMN15}.

%

\subsection{Unknown errors, robustness and interpolation}\label{ss:functions}
A drawback of the main results so far (Theorems \ref{t:nonuniformChebLeg} and \ref{t:unifquasiopt}) is that they assume the \textit{a priori} bound \R{tail_bound}, i.e.
\be{
\label{tail_bound_2}
\frac{1}{m} \sum^{m}_{j=1} \left | f(\bm{z}_j) - f_{\Lambda}(\bm{z}_j) \right |^2 \leq \eta^2,
}
for some known $\eta$.  Note that this is implied by the slightly stronger condition
\bes{
\| f - f_{\Lambda} \|_{L^\infty} \leq \eta.
}
Such an $\eta$ is required in order to formulate the optimization problem \R{weighted_l1_min} to recover $f$.  Moreover, in view of the error bounds in Theorems \ref{t:nonuniformChebLeg} and \ref{t:unifquasiopt}, one expects a poor estimation of $\eta$ to yield a larger recovery error.  Another drawback of the current approach is that the approximation $\tilde{f}$ does not interpolate $f$; a property which is sometimes desirable in applications.

We now consider the removal of the condition \R{tail_bound}.  This follows the work of \cite{BASBCScorrecting,BASBCSmodel}.  To this end, let $\eta \geq 0$ be arbitrary, i.e.\ \R{tail_bound_2} need not hold, and consider the minimization problem
\be{
\label{weighted_l1_min_new}
\min_{\bm{d} \in \bbC^n} \nm{\bm{d}}_{1,\bm{u}}\ \mbox{s.t.}\ \nm{\bm{y} - A \bm{d} }_{2} \leq \eta.
}
If $\hat{\bm{c}}_{\Lambda} = (\hat{c}_{\bm{i}} )_{\bm{i} \in \Lambda}$ is a minimzier of this problem, we define, as before, the corresponding approximation
\bes{
\tilde{f} = \sum_{\bm{i} \in \Lambda} \hat{c}_{\bm{i}} \phi_{\bm{i}}.
}
Note that if $\eta = 0$ then $\tilde{f}$ exactly interpolates $f$ at the sample points $\{\bm{z}_j \}^{m}_{j=1}$.

An immediate issue with the minimization problem \R{weighted_l1_min_new} is that the truncated vector of coefficients $\bm{c}_{\Lambda}$ is not generally feasible.  Indeed, $\bm{y} - A \bm{c}_{\Lambda} = \bm{e}_{\Lambda}$, where $\bm{e}_{\Lambda}$ is as in \R{e_Lambda_def} and is generally nonzero.  In fact, is not even guaranteed that the feasibility set of \R{weighted_l1_min_new} is nonempty.  However, this will of course be the case whenever $A$ has full rank $m$.  Under this assumption, one then has the following (see \cite{BASBCScorrecting}):

\thm{
\label{t:unifquasiopt_nobd}
Let $\epsilon$, $k$, $m$, $\gamma$, $\Lambda$, $\{ \phi_{\bm{i}} \}_{\bm{i} \in \bbN^d_0}$ and $\bm{z}_{1},\ldots,\bm{z}_m$ be as in Theorem \ref{t:unifquasiopt}.  Then with probability at least $1-\epsilon$ the following holds.  For any $\eta \geq 0$ and $f \in L^2(D) \cap L^\infty(D)$ the approximation
\bes{
\tilde{f} = \sum_{\bm{i} \in \Lambda} \hat{c}_{\bm{i}} \phi_{\bm{i}},
}
where $\bm{\hat{c}}_{\Lambda} = (\hat{c}_{\bm{i}})_{\bm{i} \in \Lambda}$ is a solution of \R{weighted_l1_min_new} with $A$ and $\bm{y}$ given by \R{y_A_def} satisfies
\be{
\label{ferr_unif_nobd}
\nm{ f - \tilde{f} }_{L^\infty} \leq \nm{ \bm{c} - \hat{\bm{c}}_{\Lambda} }_{1,\bm{u}} \lesssim \sigma_{k,L}(\bm{c})_{1,\bm{u}} + k^{\gamma/2} \left ( \eta + \nm{\bm{e}_{\Lambda}}_{2} + T_{\bm{u}}(A,\Lambda,\bm{e}_{\Lambda},\eta) \right ) 
}
and
\be{
\label{ferr_L2_nobd}
\nm{ f - \tilde{f} }_{L^2_{\nu}} = \nm{ \bm{c} - \hat{\bm{c}}_{\Lambda} }_{2}  \lesssim \frac{\sigma_{k,L}(\bm{c})_{1,\bm{u}}}{k^{\gamma/2}} + \eta + \nm{\bm{e}_{\Lambda}}_{2} + T_{\bm{u}}(A,\Lambda,\bm{e}_{\Lambda},\eta),
}
where $\bm{c} \in \ell^2(\bbN^d_0)$ are the coefficients of $f$ in the basis  $\{ \phi_{\bm{i}} \}_{\bm{i} \in \bbN^d_0}$, $\bm{e}_{\Lambda}$ is as in \R{e_Lambda_def} and
\be{
\label{tail_err}
T_{\bm{u}}(A,\Lambda,\bm{e}_{\Lambda},\eta) = \min \left \{ \frac{\nm{\bm{d}}_{1,\bm{u}}}{k^{\gamma/2}} : \bm{d} \in \bbC^n,\  \nm{A \bm{d} - \bm{e}_{\Lambda} }_{2} \leq \eta \right \}.
}
}
\prf{
We follow the steps of the proof of Theorem \ref{t:unifquasiopt} with some adjustments to take into account the fact that $\bm{c}_{\Lambda}$ may not be feasible.  First, let $S$ be such that \R{Sdef1} holds and set $\bm{d} = \bm{c}_{\Lambda} - \hat{\bm{c}}_{\Lambda}$ and $T = \Lambda \backslash S$.  Then, arguing in a similar way we see that
\eas{
\nm{\bm{d}_{T}}_{1,\bm{u}}  & \leq 2 \sigma_{k,L}(\bm{c})_{1,\bm{u}} +  \nm{ \bm{d}_{S}}_{1,\bm{u}} + \nm{\hat{\bm{c}}_{\Lambda}}_{1,\bm{u}}  - \nm{\bm{c}_{\Lambda}}_{1,\bm{u}}
\\
& \leq 2 \sigma_{k,L}(\bm{c})_{1,\bm{u}} +  \nm{ \bm{d}_{S}}_{1,\bm{u}} + \nm{\bm{g} - \bm{c}_{\Lambda}}_{1,\bm{u}},
}
where $\bm{g} \in \bbC^n$ is any point in the feasible set of \R{weighted_l1_min_new}.  By Lemma \ref{LRIPimpliesNSP} we have 
\bes{
\nm{\bm{d}_{T}}_{1,\bm{u}} \leq 2 \sigma_{k,L}(\bm{c})_{1,\bm{u}} + \rho \nm{\bm{d}_{T}}_{1,\bm{u}} + \tau \sqrt{s(k)} \nm{A \bm{d}}_{2} +  \nm{\bm{g} - \bm{c}_{\Lambda}}_{1,\bm{u}}.
}
Notice that $\nm{A \bm{d}}_{2} = \nm{\bm{y} - \bm{e}_{\Lambda} - A \hat{\bm{c}}_{\Lambda}}_{2} \leq \nm{\bm{e}_{\Lambda}}_{2} + \eta$,
and therefore
\bes{
\nm{\bm{d}_T}_{1,\bm{u}} \lesssim \sigma_{k,L}(\bm{c})_{1,\bm{u}} + \sqrt{s(k)} \left ( \nm{\bm{e}_{\Lambda}}_{2} + \eta \right ) + \nm{\bm{g} - \bm{c}_{\Lambda}}_{1,\bm{u}}.
}
Hence, by similar arguments, it follows that
\be{
\label{l1u_bound}
\nm{\bm{c} - \bm{c}_{\Lambda}}_{1,\bm{u}} \lesssim \sigma_{k,L}(\bm{c})_{1,\bm{u}} + k^{\gamma/2} \left ( \nm{\bm{e}_{\Lambda}}_{2} + \eta \right ) + \nm{\bm{g} - \bm{c}_{\Lambda}}_{1,\bm{u}},
}
for any feasible point $\bm{g}$.  After analogous arguments, we also deduce the following bound in the $\ell^2$-norm:
\be{
\label{l2_bound}
\nm{\bm{c} - \bm{c}_{\Lambda}}_2 \lesssim \frac{\sigma_{k,L}(\bm{c})_{1,\bm{u}}}{k^{\gamma/2}} + \left ( \nm{\bm{e}_{\Lambda}}_{2} + \eta \right ) + k^{-\gamma/2} \nm{\bm{g} - \bm{c}_{\Lambda}}_{1,\bm{u}}.
}
To complete the proof, we consider the term $\nm{\bm{g} - \bm{c}_{\Lambda}}_{1,\bm{u}}$.  Write $\bm{g} = \bm{c}_{\Lambda} + \bm{g}'$ and notice that $\bm{g}$ is feasible if and only if $\bm{g}'$ satisfies $\nm{A \bm{g}' - \bm{e}_{\Lambda}} \leq \eta$.  Since $\bm{g}'$ is arbitrary we get the result.
\qed
}

The two error bounds \R{ferr_unif_nobd} and \R{ferr_L2_nobd} in this theorem are analogous to \R{ferr_unif} and \R{ferr_L2} in Theorem \ref{t:unifquasiopt}.  They remove the condition $\eta \geq \nm{\bm{e}_{\Lambda}}_{2}$ at the expense of an additional term $T_{\bm{u}}(A,\Lambda,\bm{e}_{\Lambda},\eta)$.  We now provide a bound for this term (see \cite{BASBCScorrecting}):

\thm{
Consider the setup of Theorem \ref{t:unifquasiopt_nobd} and let $T_{\bm{u}}(A,\Lambda,\bm{e}_{\Lambda},\eta)$ be as in \R{tail_err}.  If $A$ has full rank, then
\be{
\label{Tbound}
T_{\bm{u}}(A,\Lambda,\bm{e}_{\Lambda},\eta) \leq  \frac{k^{\alpha/2} \sqrt{L}}{\sigma_{\min} \left ( \sqrt{\frac{m}{n}} A^* \right ) }\max \left \{ \nm{\bm{e}_{\Lambda}}_2 - \eta , 0 \right \},
}
where $L$ is as in \R{uniformlogfactor} and $\alpha = 1,2$ in the Chebyshev or Legendre cases respectively.
}
\prf{
If $\eta \geq \nm{\bm{e}_{\Lambda}}_2$ then the result holds trivially.  Suppose now that $\eta < \nm{\bm{e}_{\Lambda}}_2$.  Since $\nm{\bm{e}_{\Lambda}}_2 \neq 0$ in this case, we can define $\bm{d} = (1 - \eta / \nm{\bm{e}_{\Lambda}}_2 ) A^{\dag} \bm{e}_{\Lambda}$, where $A^{\dag}$ denotes the pseudoinverse.  Then $\bm{d}$ satisfies $\nm{A \bm{d} - \bm{e}_{\Lambda}}_{2} = \eta$, and therefore
\bes{
k^{\gamma/2} T_{\bm{u}}(A,\Lambda,\bm{e}_{\Lambda},\eta) \leq \nm{d}_{1,\bm{u}} \leq \sqrt{| \Lambda |_{\bm{u}}} \nm{\bm{d}}_{2} \leq \frac{\sqrt{| \Lambda |_{\bm{u}}}}{\sigma_{\min}(A^*)} \left ( \nm{\bm{e}_{\Lambda}}_2- \eta  \right ).
}
Equation \R{uniformsamplecomplexity} implies that $\sqrt{\frac{m}{k^{\gamma}}} \lesssim \sqrt{L}$, and hence
\be{
\label{Tbound1}
T_{\bm{u}}(A,\Lambda,\bm{e}_{\Lambda},\eta) \lesssim \sqrt{\frac{| \Lambda |_{1,\bm{u}}}{n} } \frac{\sqrt{L}}{\sigma_{\min}  \left ( \sqrt{\frac{m}{n}} A^* \right )} \left (\nm{\bm{e}_{\Lambda}}_2- \eta  \right ).
}
It remains to estimate $| \Lambda |_{1,\bm{u}}$.  For the Chebyshev case, we apply \R{Chebu} to give
\bes{
| \Lambda |_{1,\bm{u}} = \sum_{\bm{i} \in \Lambda} 2^{\nm{\bm{i}}_0} \leq \sum_{\bm{i} \in \Lambda} \prod^{d}_{j=1} \left ( i_j + 1 \right ) \leq k \sum_{\bm{i} \in \Lambda} 1  = k n
}
where in the penultimate step we used the definition of the hyperbolic cross \R{HCdef}.  For the Legendre case, we use \R{Legu} to get
\bes{
| \Lambda |_{1,\bm{u}} = \sum_{\bm{i} \in \Lambda} \prod^{d}_{j=1} \left ( 2 i_j + 1 \right ) \leq \sum_{\bm{i} \in \Lambda} 2^{\nm{\bm{i}}_0} \prod^{d}_{j=1} \left ( i_j + 1 \right ) \leq k^2 n.
}
This completes the proof. 
\qed
}

The error bound \R{Tbound} suggests that the effect of removing the condition $\eta \geq \nm{\bm{e}_{\Lambda}}_{2}$ is at most a small algebraic factor in $k$, a log factor and term depending on the minimal singular value of the scaled matrix $\sqrt{\frac{m}{n}} A^*$.  We discuss this latter term further in below.  Interestingly, this bound suggests that a good estimate of $\nm{\bm{e}_{\Lambda}}_2$ (when available) can reduce this error term.  Indeed, one has $T_{\bm{u}}(A,\Lambda,\bm{e}_{\Lambda},\eta) \rightarrow 0$ linearly in $\nm{\bm{e}_{\Lambda}}_2-\eta \rightarrow 0^{+}$.  Hence estimation procedures aiming to tune $\eta$  -- for example, cross validation (see \S \ref{ss:numerical}) -- are expected to yield reduced error over the case $\eta = 0$, for example.

It is beyond the scope of this chapter to provide theoretical bounds on the minimal singular value of the scaled matrix  $\sqrt{\frac{m}{n}} A^*$.  We refer to \cite{BASBCSmodel} for a more comprehensive treatment of such bounds.  However, we note that it is reasonable to expect that $\sigma_{\min}(\sqrt{\frac{m}{n}} A^*) \approx 1$ under appropriate conditions on $m$ and $n$.  Indeed:

\lem{
Let $B = \bbE \left ( \frac{m}{n} A A^* \right )$, where $A$ is the matrix of Theorem \ref{t:unifquasiopt_nobd}.  Then the minimal eigenvalue of $B$ is precisely $1-1/n$.
}
\prf{
We have
$
\bbE \left ( \frac{m}{n} A A^* \right )_{j,l} = \bbE \left ( \frac{1}{n} \sum_{\bm{i} \in \Lambda} \phi_{\bm{i}}(\bm{z}_j) \phi_{\bm{i}}(\bm{z}_l) \right ).
$
When $l = j$ this gives $\bbE \left ( \frac{m}{n} A A^* \right )_{j,j} = 1$.
Conversely, since $\{ \phi_{\bm{i}} \}_{\bm{i} \in \bbN^d_0}$ are orthogonal polynomials one has $\int_{D} \phi_{\bm{i}}(\bm{z}) \D \nu = \ip{\phi_{\bm{i}} }{\phi_{\bm{0}} }_{L^2_{\nu}} = \delta_{\bm{i},\bm{0}}$, and therefore for $l \neq j$ one has
$
\bbE \left ( \frac{m}{n} A A^* \right )_{j,l} = \frac{1}{n} \sum_{\bm{i} \in \Lambda} \left ( \int_{D} \phi_{\bm{i}}(\bm{z}) \D \nu \right )^2  = \frac{1}{n},
$
It is now a straightforward calculation to show that $\lambda_{\min}(B) = 1-1/n$.
\qed
}

\rem{
\label{rem:computable}
Although complete theoretical estimates $T_{\bm{u}}(A,\Lambda,\bm{e}_{\Lambda},\eta)$ are outside the scope of this work, it is straightforward to derive a bound that can be computed.  Indeed, it follows immediately from \R{Tbound1} that
\bes{
T_{\bm{u}}(A,\Lambda,\bm{e}_{\Lambda},\eta) \lesssim Q_{\bm{u}}(A) \sqrt{L} \max \left \{  \nm{\bm{e}_{\Lambda}}_2- \eta  , 0 \right \},
}
where
\be{
\label{QuA}
Q_{\bm{u}}(A) = \sqrt{\frac{| \Lambda |_{1,\bm{u}}}{n}} \frac{1}{\sigma_{\min}\left ( \sqrt{\frac{m}{n}} A^* \right )}.
}
Hence, up to the log factor, the expected robustness of \R{weighted_l1_min_new} can be easily checked numerically.  See \S \ref{ss:numerical} for some examples of this approach.
}

\rem{\label{rem:err_sources} For pedagogical reasons, we have assumed the truncation of $f$ to  $f_\Lambda$ is the only source of error $\bm{e}_\Lambda$ affecting the measurements $\bm{y}$ (recall \eqref{e_Lambda_def}). There is no reason for this to be the case, and $\bm{e}_{\Lambda}$ may incorporate other errors without changing any of the above results.  We note that concrete applications often give rise to other sources of unknown error.  For example, in UQ, we usually aim at approximating a function of the form $f(\bm{z}) = q(u(\bm{z}))$, where  $u(\bm{z})$ is the solution to a PDE depending on some random coefficients $\bm{z}$ and $q$ is a quantity of interest (see \cite{DoostanOwhadiSparse,KarniadakisUQCS}, for example). In this case, each sample $f(\bm{z}_j)$ is typically subject to further sources of inaccuracy, such as the numerical error associated with the PDE solver employed to compute $u(\bm{z}_j)$ (e.g.\ a finite element method) and, possibly, the error committed  evaluating $q$ on $u(\bm{z}_j)$ (e.g.\ numerical integration).
}

\rem{Our analysis based on the estimation of the tail error \eqref{tail_err} can be compared with the robustness analysis of basis pursuit based on the so-called \emph{quotient property} \cite{FoucartRauhutCSbook}. However, this analysis is limited to the case of \emph{basis pursuit}, corresponding to the optimization program \eqref{weighted_l1_min_new} with $\bm{u} = \bm{1}$ (i.e., unweighted $\ell^1$ norm) and $\eta = 0$. In the context of compressed sensing, random matrices that are known to fulfill the quotient property with high probability are gaussian, subgaussian, and Weibull matrices \cite{Foucart2014,Wojtaszczyk2010}.  For further details we refer to \cite{BASBCSmodel}.}

\subsection{Numerical results}
\label{ss:numerical}

We conclude this chapter with a series of numerical results.  First, in Figure \ref{f:ErrLU} and \ref{f:ErrCC} we show the approximation of several functions via weighted $\ell^1$ minimization.  Weights of the form $w_{\bm{i}} = (u_{\bm{i}})^{\alpha}$ are used for several different choices of $\alpha$.  In agreement with the discussion in \S \ref{ss:nonuniform} the choice $\alpha = 1$, i.e.\ $w_{\bm{i}} = u_{\bm{i}}$ generally gives amongst the smallest error.  Moreover, while larger values of $\alpha$ sometime give a smaller error, this is not the case for all functions.  Notice that in all cases unweighted $\ell^1$ minimization gives a worse error than weighted $\ell^1$ minimization.  As is to be expected, the improvement offered by weighted $\ell^1$ minimization in the Chebyshev case is less significant in moderate dimensions than for Legendre polynomials.

\begin{figure}
\begin{center}
\begin{tabular}{ccc}
\includegraphics[width=5.25cm]{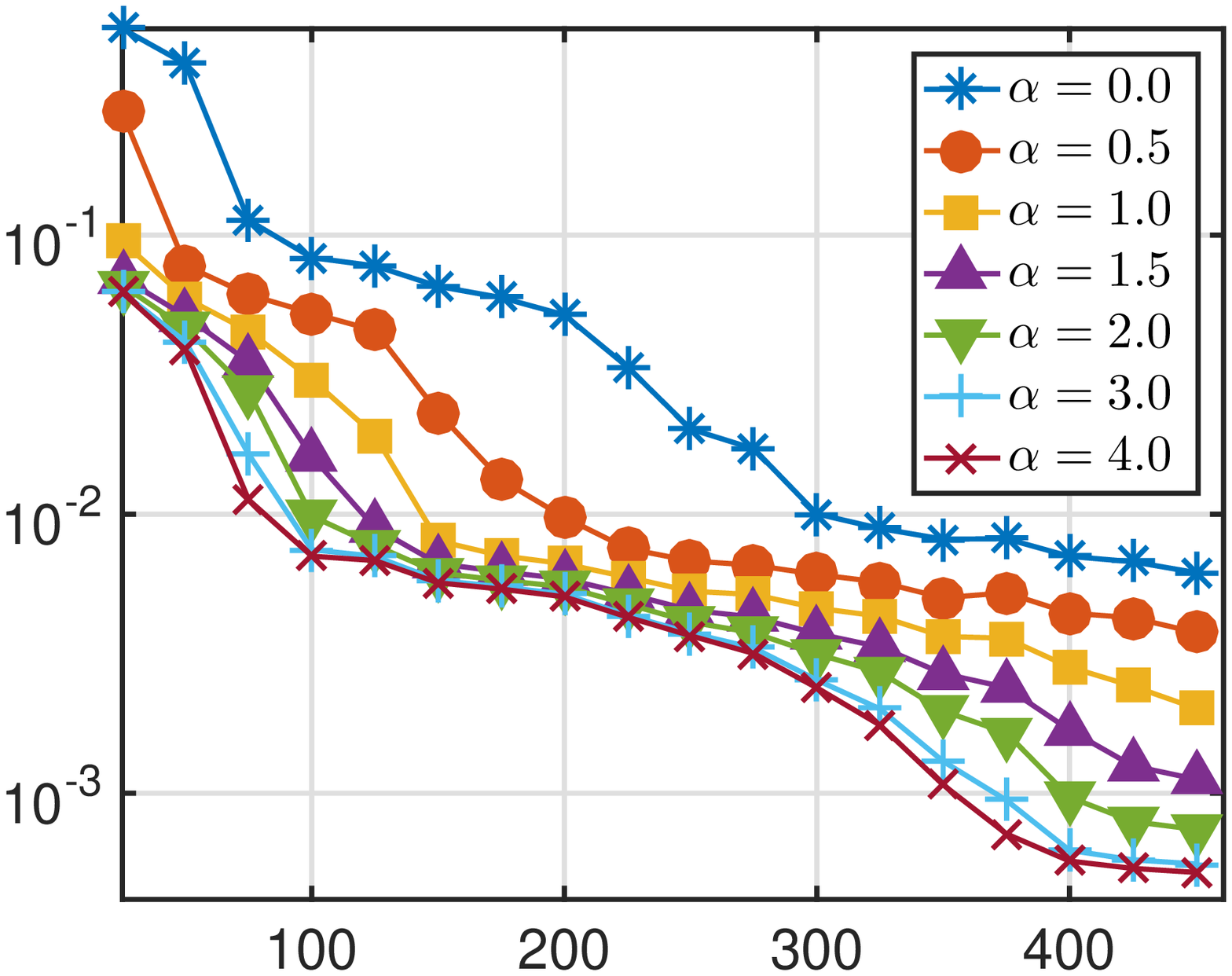} &&
\includegraphics[width=5.25cm]{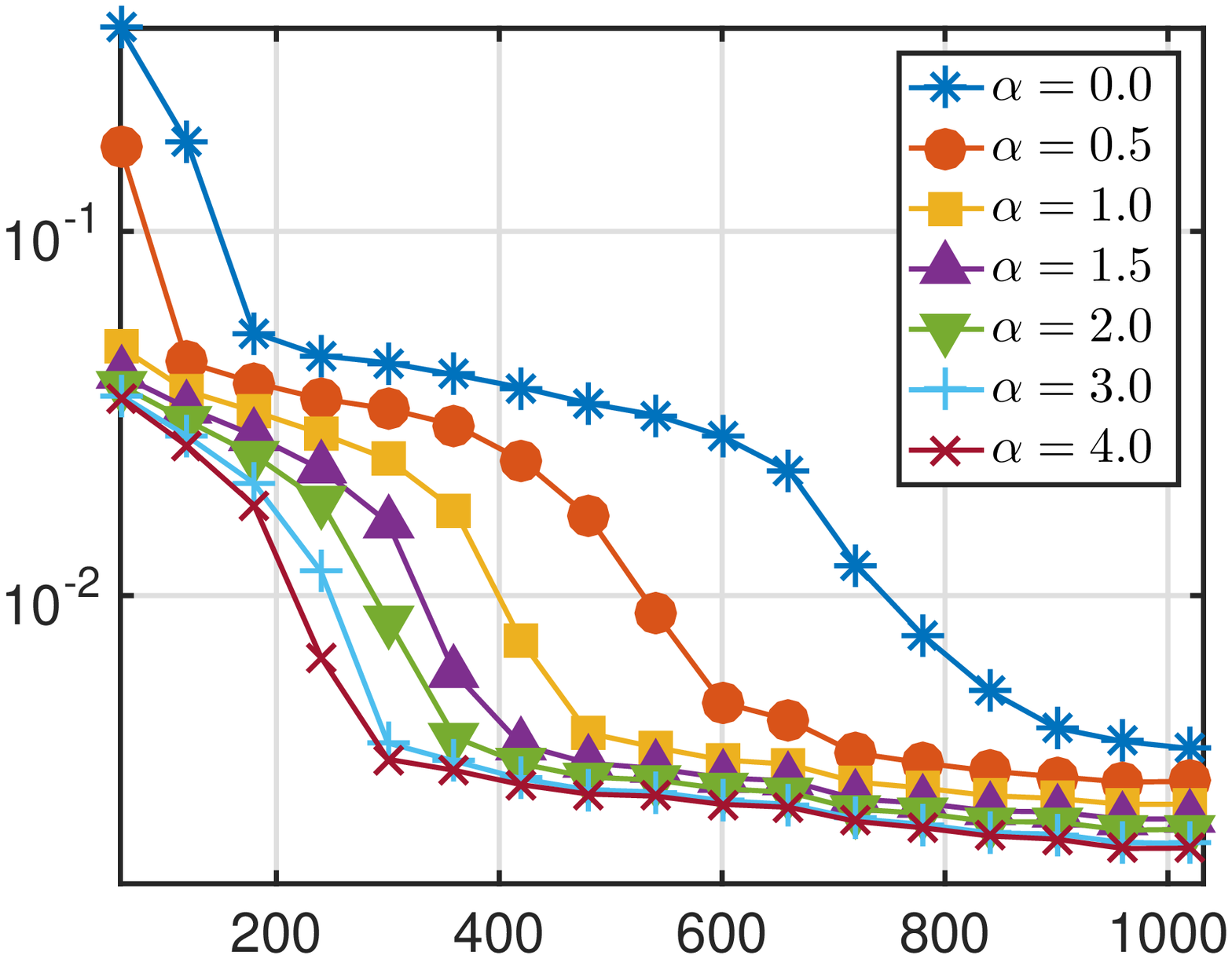} 
\\
\includegraphics[width=5.25cm]{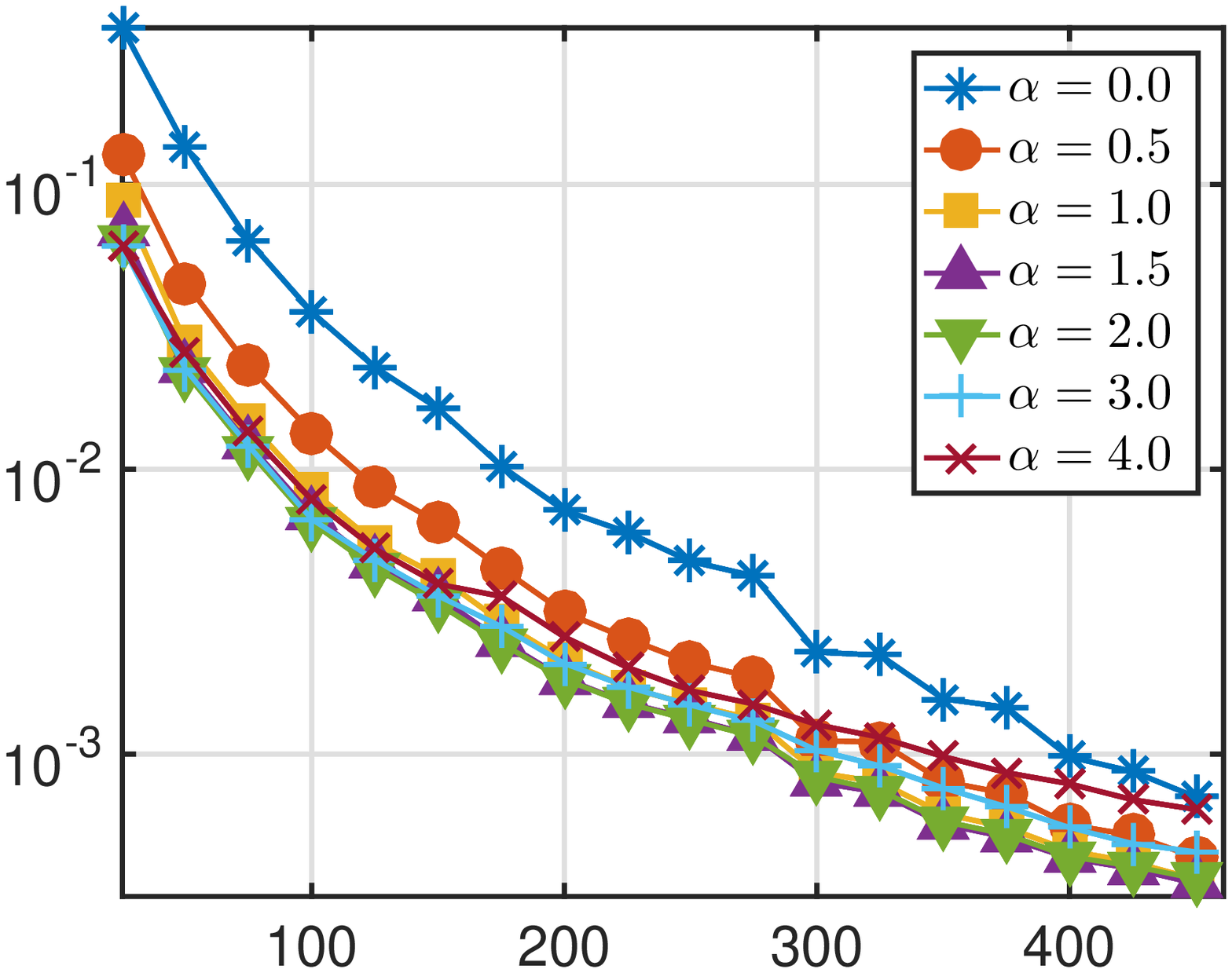} &&
\includegraphics[width=5.25cm]{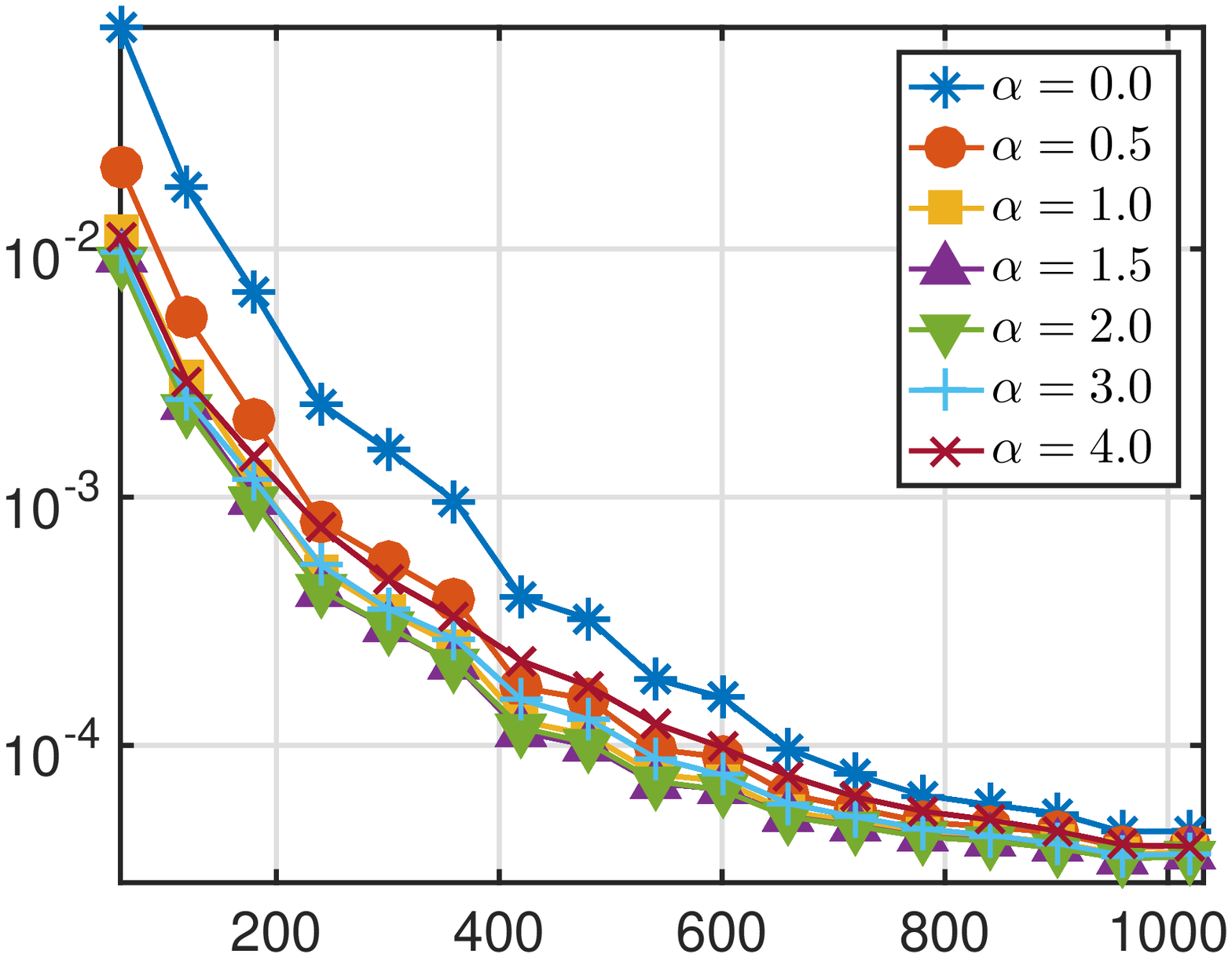} 
\\
$(d,k,n) = (8,22,1843)$ && $(d,k,n) = (16,13,4129)$
\end{tabular}
\end{center}
\caption{
The error $\| f - \tilde{f} \|_{L^\infty}$ (averaged over $50$ trials) against $m$ for Legendre polynomials.  Here $\tilde{f} = \sum_{\bm{i} \in \Lambda} \hat{c}_{\bm{i}} \phi_{\bm{i}}$, where $\hat{\bm{c}}_{\Lambda}$ is a solution of \R{weighted_l1_min} with weights $w_{\bm{i}} = (u_{\bm{i}})^{\alpha}$ and $\Lambda = \Lambda^{\mathrm{HC}}_{k}$ a hyperbolic cross index set.  The functions used were $f(\bm{y}) = \prod^{d}_{k=d/2+1} \cos(16 y_k / 2^k) / \prod^{d/2}_{k=1} (1-y_k/4^k )$ and $f(\bm{y}) = \exp \left ( -\sum^{d}_{k=1} y_k / (2d) \right )$ (top and bottom respectively).  The weighted $\ell^1$ minimization problem was solved using the SPGL1 package \cite{spgl1:2007,BergFriedlander:2008} with a maximum of 100,000 iterations and $\eta = 10^{-12}$.
}
\label{f:ErrLU}
\end{figure}

\begin{figure}
\begin{center}
\begin{tabular}{ccc}
\includegraphics[width=5.25cm]{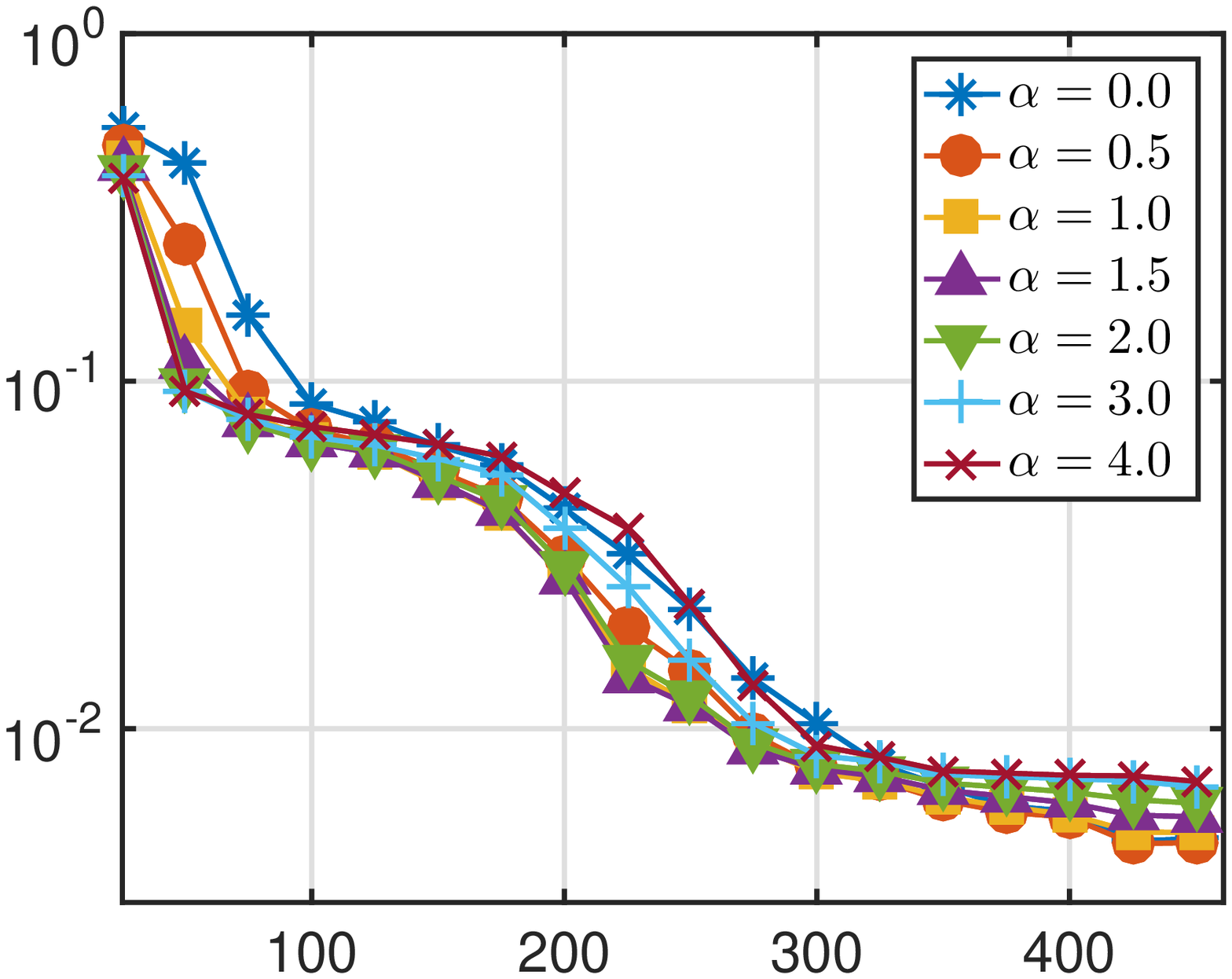} &&
\includegraphics[width=5.25cm]{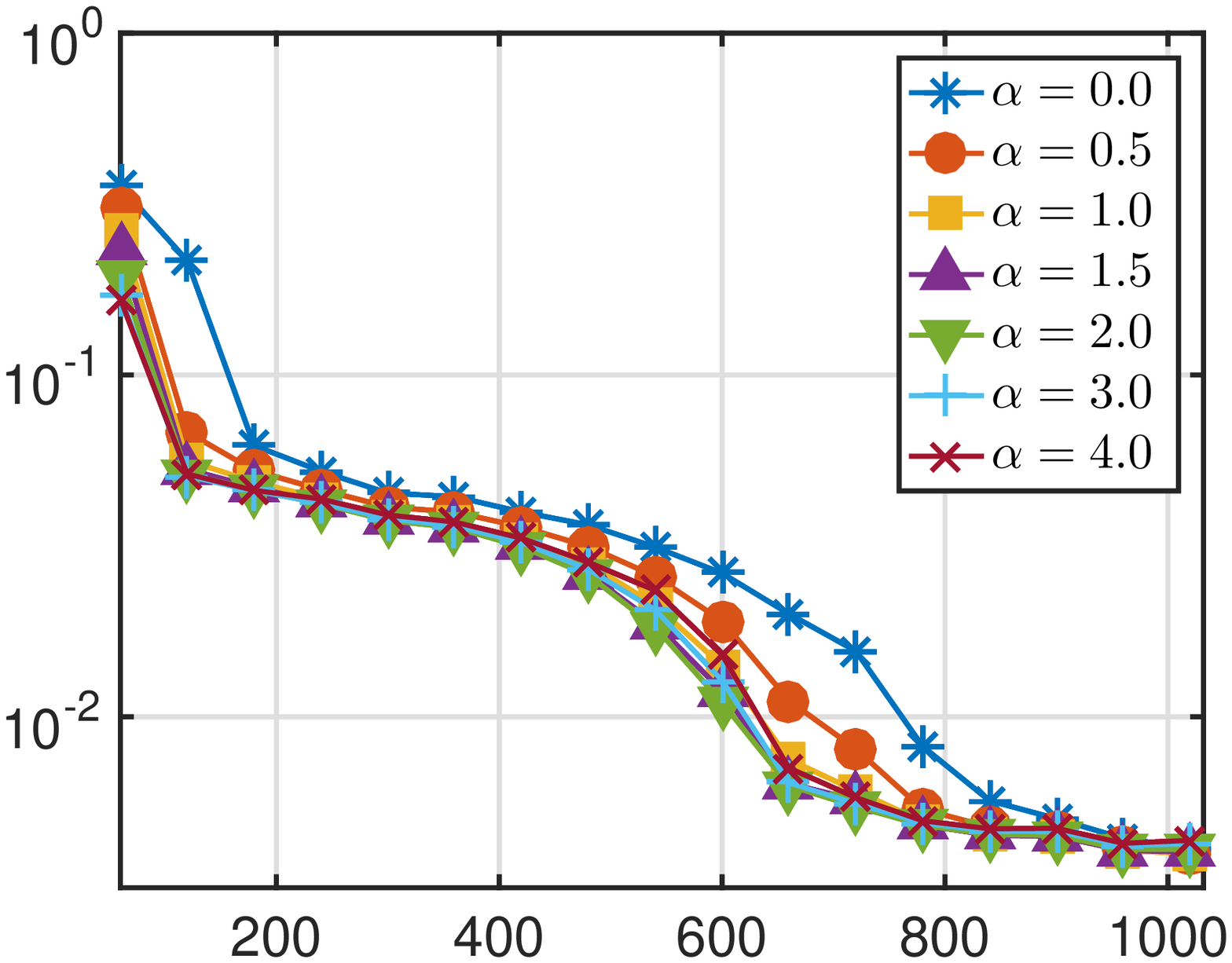} 
\\
\includegraphics[width=5.25cm]{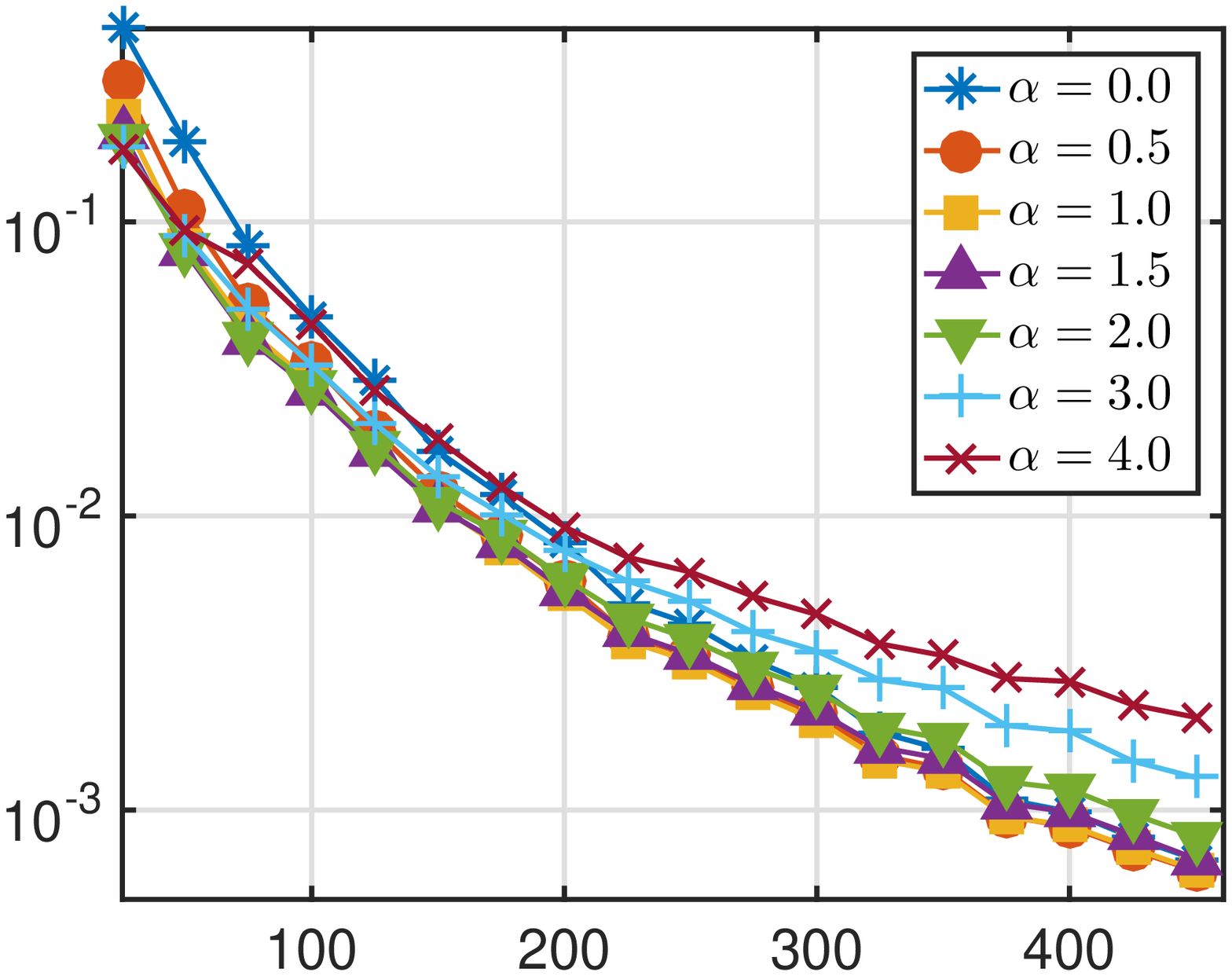} &&
\includegraphics[width=5.25cm]{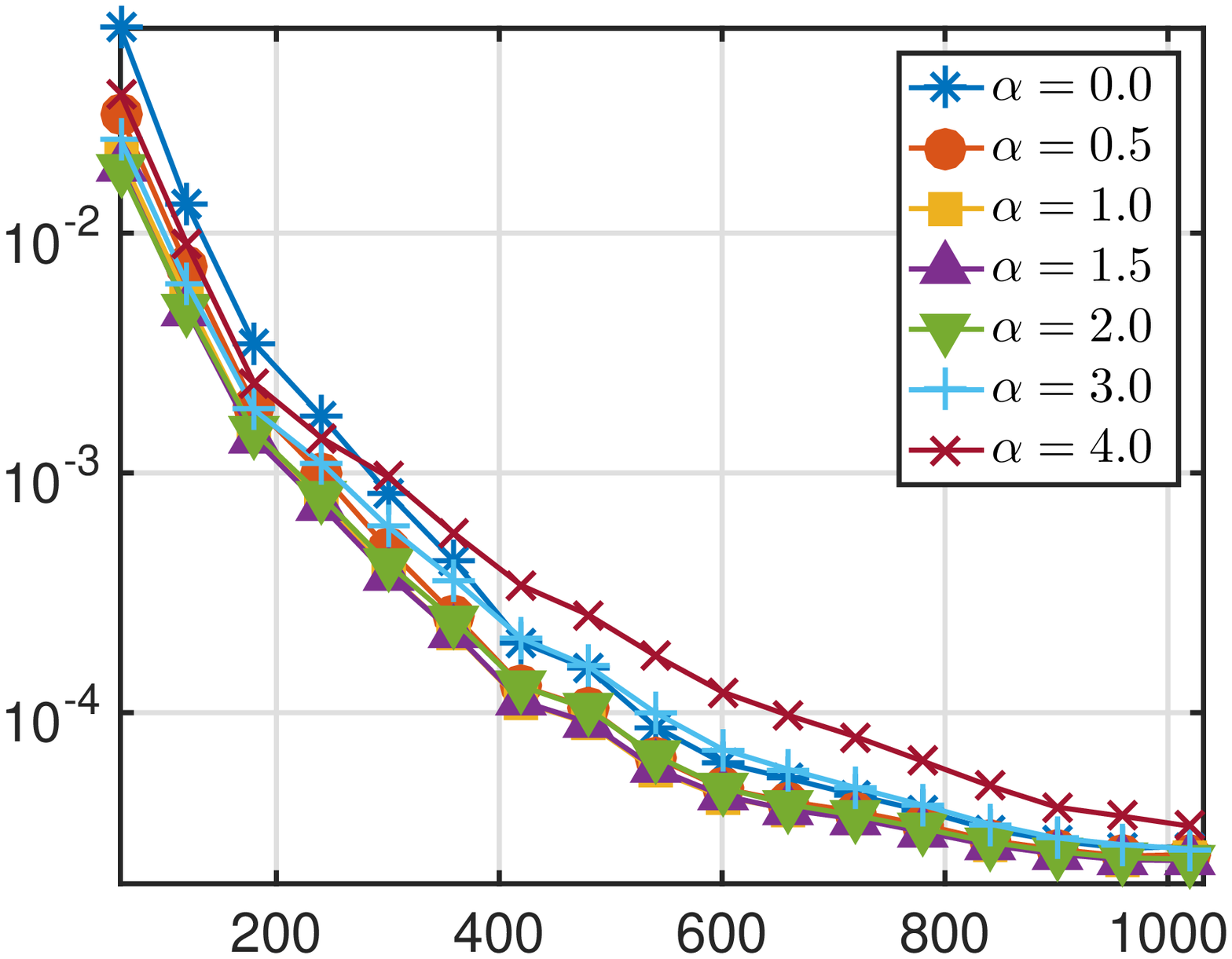} 
\\
$(d,k,n) = (8,22,1843)$ && $(d,k,n) = (16,13,4129)$
\end{tabular}
\end{center}
\caption{
The same as Figure \ref{f:ErrLU} but with Chebyshev polynomials. 
}
\label{f:ErrCC}
\end{figure}

The results in Figures \ref{f:ErrLU} and \ref{f:ErrCC} were computed by solving weighted $\ell^1$ minimization problems with $\eta$ set arbitrarily to $\eta = 10^{-12}$ (we make this choice rather than $\eta = 0$ to avoid potential infeasibility issues in the solver).  In particular, the condition \R{tail_bound_2} is not generally satisfied.  Following Remark \ref{rem:computable}, we next assess the size of the constant $Q_{\bm{u}}(A)$ defined in \R{QuA}.  Table \ref{t:ErrConst} shows the magnitude of this constant for the setups considered in Figures \ref{f:ErrLU} and \ref{f:ErrCC}.  Over all ranges of $m$ considered, this constant is never more than 20 in magnitude.  That is to say, the additional effect due to the unknown truncation error $\nm{\bm{e}_{\Lambda}}_2$ is relatively small.

\begin{table}
\centering
\begin{tabular}{|l|c|cccccccc|}
\hline
 &   $m$   &   125 & 250 & 375 & 500 & 625 & 750 & 875 & 1000   \\
\cline{2-10}
\multirow{2}{*}{$(d,k,n) = (8,22,1843)$} & Chebyshev & 2.65 & 3.07 & 3.53 & 3.95 & 4.46 & 5.03 & 5.78 & 6.82  \\ 
& Legendre &  6.45 & 7.97 & 8.99 & 10.5 & 12.1 & 13.7 & 15.8 & 18.6 \\ 
\hline
&   $m$   &   250 & 500 & 750 & 1000 & 1250 & 1500 & 1750 & 2000 \\
\cline{2-10}
\multirow{2}{*}{$(d,k,n) = (16,13,4129)$} & Chebyshev &  2.64 & 2.93 & 3.30 & 3.63 & 3.99 & 4.41 & 4.95 & 5.62   \\ 
& Legendre & 5.64 & 6.20 & 6.85 & 7.60 & 8.32 & 8.99 & 10.1 & 11.1 \\ 
\hline
\end{tabular}
\label{t:ErrConst}
\caption{The constant $Q_{\bm{u}}(A)$ (averaged over 50 trials) for the setup considered in Figures \ref{f:ErrLU} and \ref{f:ErrCC}.}
\end{table}

In view of Remark~\ref{rem:err_sources}, in Figure~\ref{f:m_vs_err_noise} we assess the performance of weighted $\ell^1$ minimization in the presence of external sources of error corrupting the measurements. In order to model this scenario, we consider the problem \R{weighted_l1_min} where the vector of measurements is corrupted by additive noise
\begin{equation}
\label{ynoisy}
\bm{y} = \frac{1}{\sqrt{m}} (f(\bm{z}_j))_{j=1}^m + \bm{n},
\end{equation}
or, equivalently, by recalling \R{underdet_system}, 
\begin{equation}
\bm{y} = A \bm{c}_{\Lambda} + \bm{e}_\Lambda + \bm{n}.
\end{equation}
We randomly generate the noise as $\bm{n} = 10^{-3}\bm{g}/\|\bm{g}\|_2$, where $\bm{g} \in \bbR^m$ is a standard random gaussian vector, so that  $\|\bm{n}\|_2 = 10^{-3}$. Considering weights $\bm{w} = (u_{\bm{i}}^\alpha)_{\bm{i} \in \Lambda}$, with  $\alpha=0,1$, we compare the error obtained when the parameter $\eta$ in \R{weighted_l1_min} is chosen according to each of the following three strategies:
\begin{enumerate}
\item $\eta = 0$, corresponding to enforcing the exact constraint $A \bm{d} = \bm{y}$ in \R{weighted_l1_min};
\item $\eta = \eta_{oracle} = \|A \hat{\bm{c}}_{oracle} - \bm{y}\|_2$, where  $f_{oracle} = \sum_{\bm{i} \in \Lambda} (\hat{c}_{oracle})_{\bm{i}} \phi_{\bm{i}}$ is the oracle least-squares solution based on $10 n$ random samples of $f$ distributed according to $\nu$;
\item $\eta$ is estimate using a cross validation approach, as described in \cite[Section 3.5]{DoostanOwhadiSparse}, where the search of $\eta$ is restricted to the values of the form $10^{k} \cdot \eta_{oracle}$, where $k$ belongs to a uniform grid of $11$ equispaced points on the interval $[-3,3]$, $3/4$ of the samples are used as reconstruction samples and $1/4$ as validation samples. 
\end{enumerate}
The results are in accordance with the estimate \R{ferr_L2_nobd}. Indeed, as expected, for any value of $\alpha$, the recovery error associated with $\bm{n} = \bm{0}$ and $\eta = 0$ is always lower than the recovery error associated with $\bm{n}\neq \bm{0}$ and any choice of $\eta$. This can be explained by the fact that, in the right-hand side of \R{ferr_L2_nobd}, the terms $\sigma_{k,L}(\bm{c})/k^{\gamma/2}$ and $\|\bm{e}_\Lambda\|_2$ are dominated by $\eta + T_{\bm{u}}(A,\Lambda,\bm{e}_\Lambda,\eta)$ when $\bm{n}\neq \bm{0}$. Moreover, estimating $\eta$ via oracle least-squares (strategy 2) gives better results than cross validation (strategy 3), which in turn is better than the neutral choice $\eta =0$ (strategy 1). Finally, we note that the discrepancy among the three strategies is accentuated as $\alpha$ gets larger.

\begin{figure}
\begin{center}
\begin{tabular}{ccc}
\includegraphics[width=5.25cm]{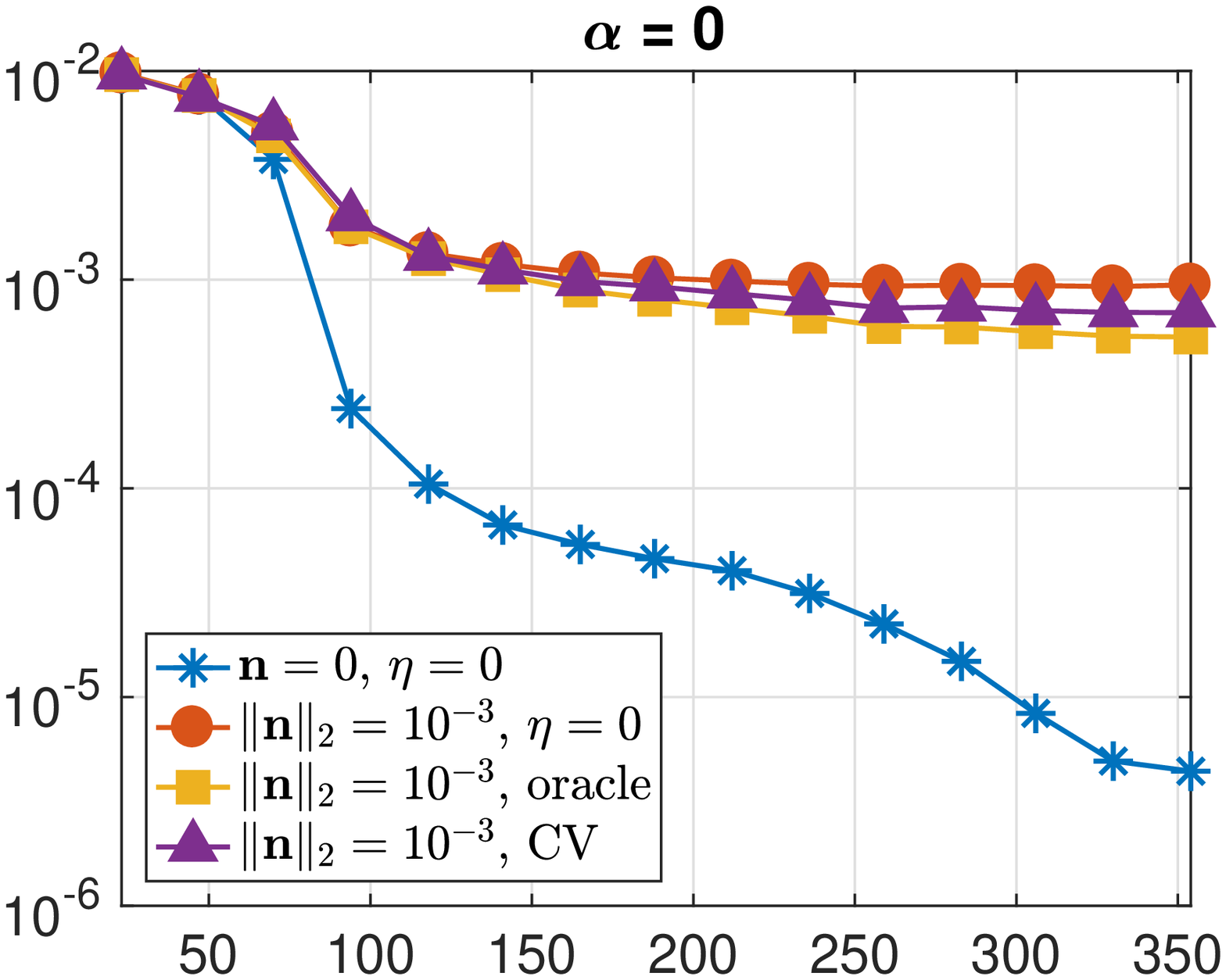} &&
\includegraphics[width=5.25cm]{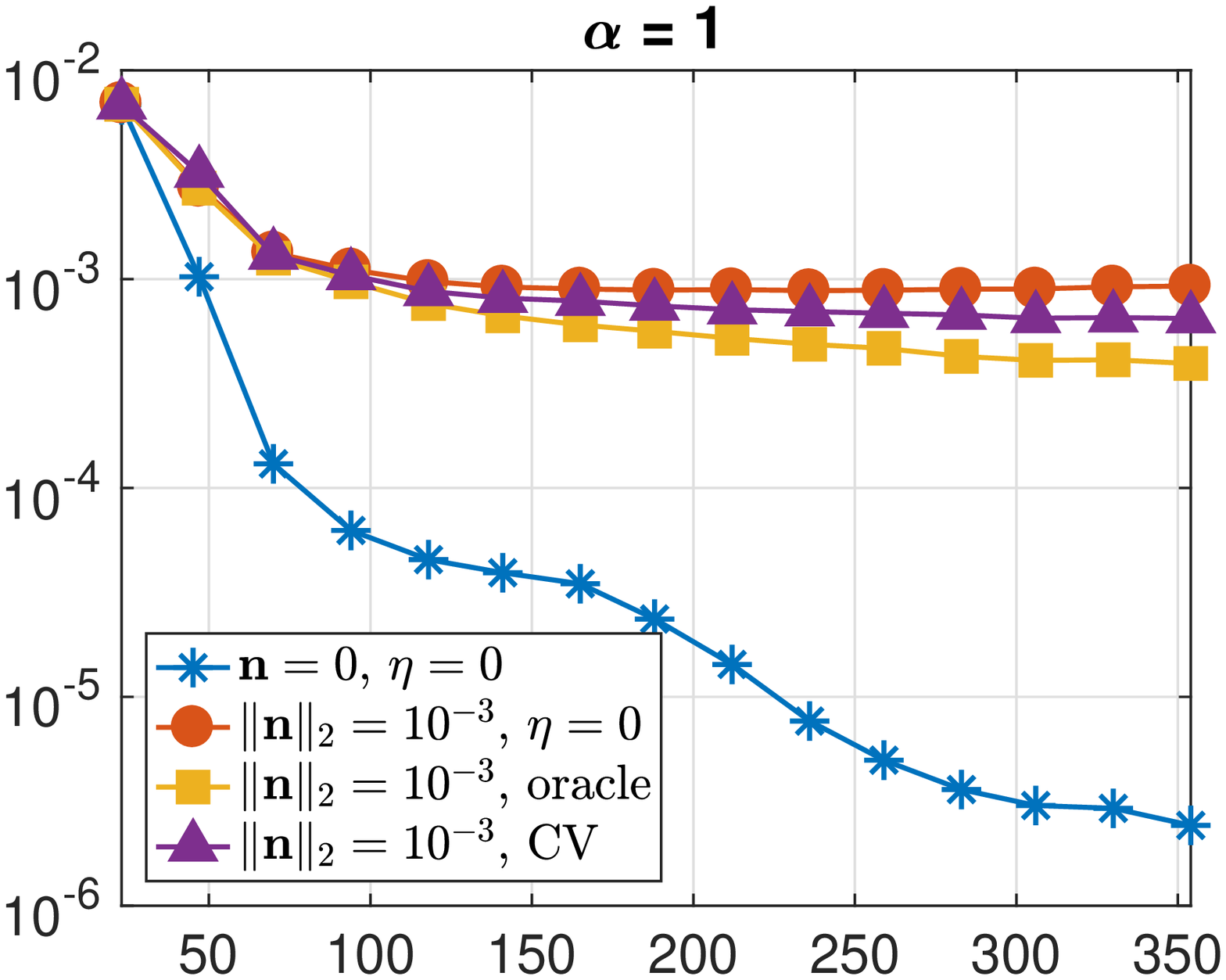}\\
\includegraphics[width=5.25cm]{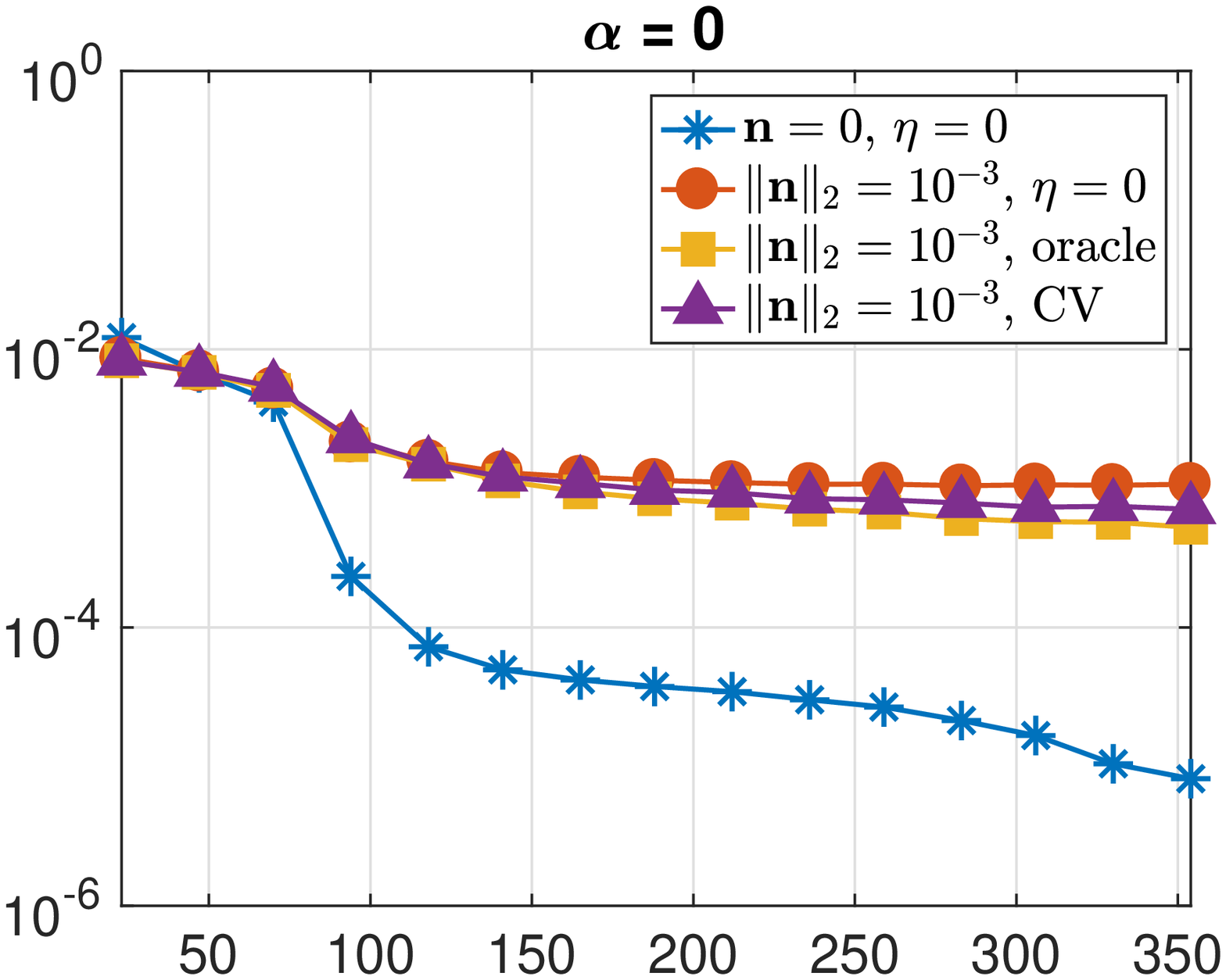} &&
\includegraphics[width=5.25cm]{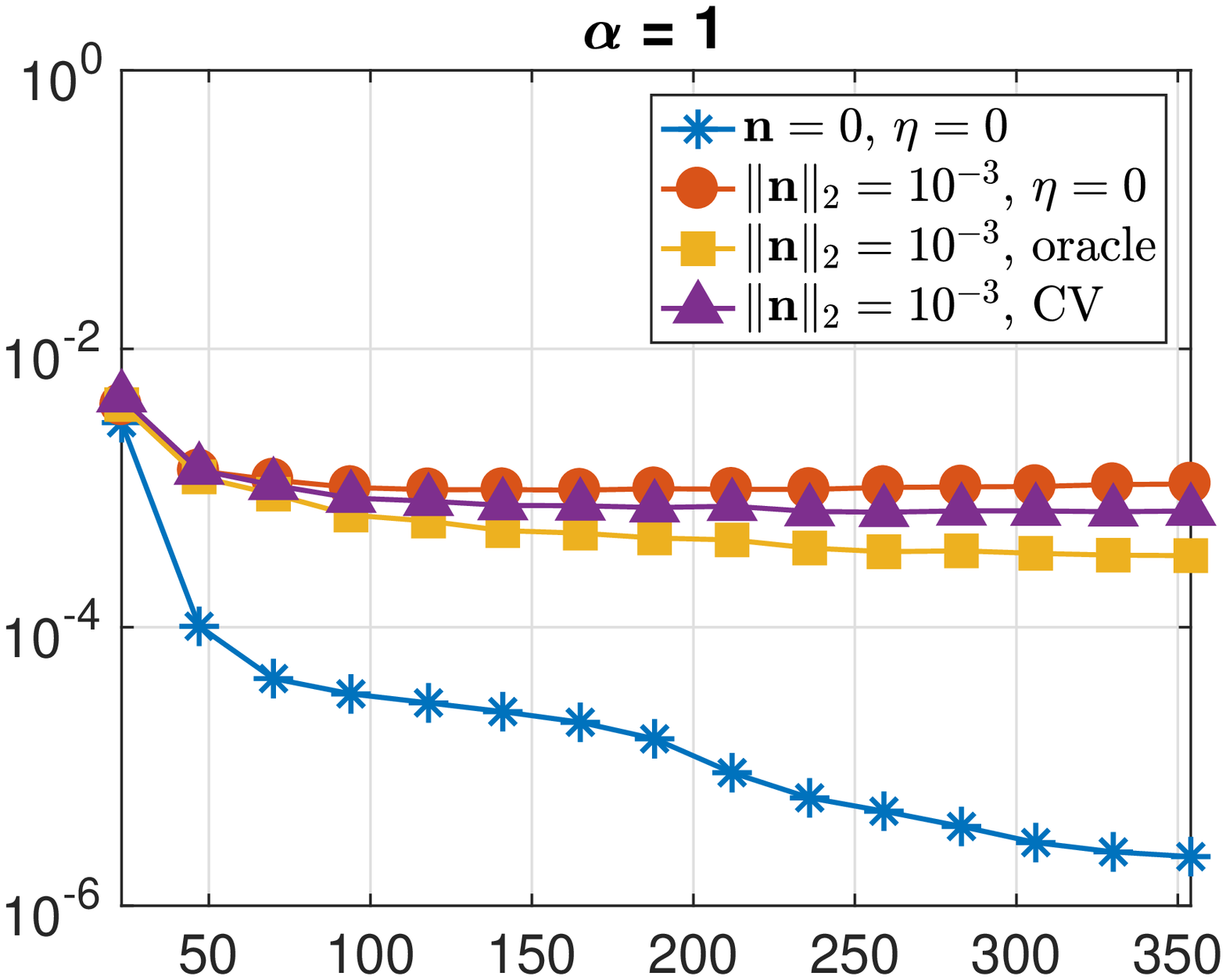}
\end{tabular}
\end{center}
\caption{\label{f:m_vs_err_noise}The error $\| f - \tilde{f} \|_{L^2_\nu}$ against $m$.  Here $\tilde{f} = \sum_{\bm{i} \in \Lambda} \hat{c}_{\bm{i}} \phi_{\bm{i}}$, where $\hat{\bm{c}}_{\Lambda}$ is a solution of \R{weighted_l1_min} with weights $\bm{w} = (u^{\alpha}_{\bm{i}})_{\bm{i}\in\Lambda}$, with $\alpha = 0$ (left) and $\alpha = 1$ (right), and $\bm{y}$ defined as in \eqref{ynoisy}. Regarding $\{\phi_{\bm{i}}\}_{\bm{i}\in\Lambda}$ and $\nu$, the Chebyshev polynomials with the Chebyshev measure are employed in the top line and the Legendre polynomials with the uniform measure in the bottom line.
We choose $d = 8$ and $\Lambda = \Lambda^{\text{HC}}_{19}$ with $n = |\Lambda| = 1771$. For each value of $m$, we average the error over 50 trials considering three different strategies for the choice of $\eta$: namely, $\eta = 0$, estimation via oracle least-squares, and cross validation (CV). The function approximated is $f(\bm{y}) = \exp\left ( -\sum^{d}_{k=1} \cos(y_k) / (8 d) \right )$. 
}
\end{figure}

In the next experiment we highlight the importance of the parameter $\eta$ when solving \R{weighted_l1_min} with measurements subject to external sources of error (recall Remark~\ref{rem:err_sources}). We corrupt the measurements by adding random noise $\bm{n}$ with norm $\|\bm{n}\|_2=10^{-3}$, analogously to \R{ynoisy}. Then, for different values of $\eta$ from $10^{-5}$ to $10$ we solve \R{weighted_l1_min} with weights $\bm{w}=(\bm{u}_{\bm{i}}^\alpha)_{\bm{i}\in\Lambda}$ and $\alpha = 0,1$. The resulting recovery errors with respect to the $L^2_\nu$ norm (averaged over 50 trials) are plotted as a function of $\eta$ in Figure~\ref{f:eta_vs_err_noise}. For every value of $\alpha$, the resulting curve is constant for the smallest and largest values of $\eta$. In between, the curve exhibits a global minimum, which corresponds to an optimal calibration of $\eta$. The values of $\eta$ estimated via oracle least-squares and cross validation are both able to approximate the global minimum on average. However, cross validation has a larger  standard deviation compared to the former (see Table~\ref{t:mean_std_eta}). This explains why the performance of cross validation is suboptimal in Figure~\ref{f:m_vs_err_noise}. We also notice that the global minimum is more pronounced as $\alpha$ gets larger, in accordance to the observations in Figure~\ref{f:m_vs_err_noise}.
\begin{figure}
\begin{center}
\begin{tabular}{ccc}
\includegraphics[width=5.25cm]{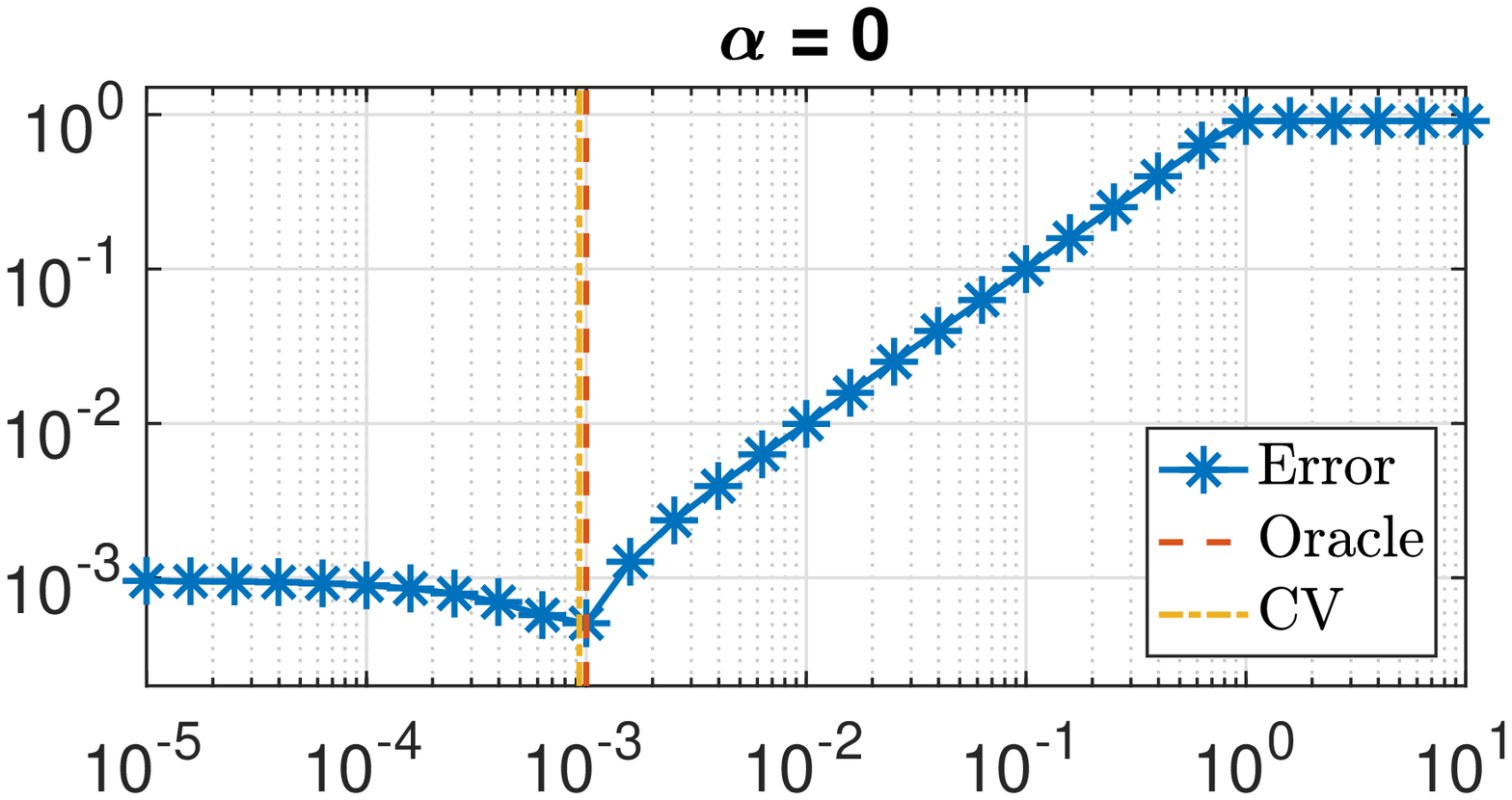} && 
\includegraphics[width=5.25cm]{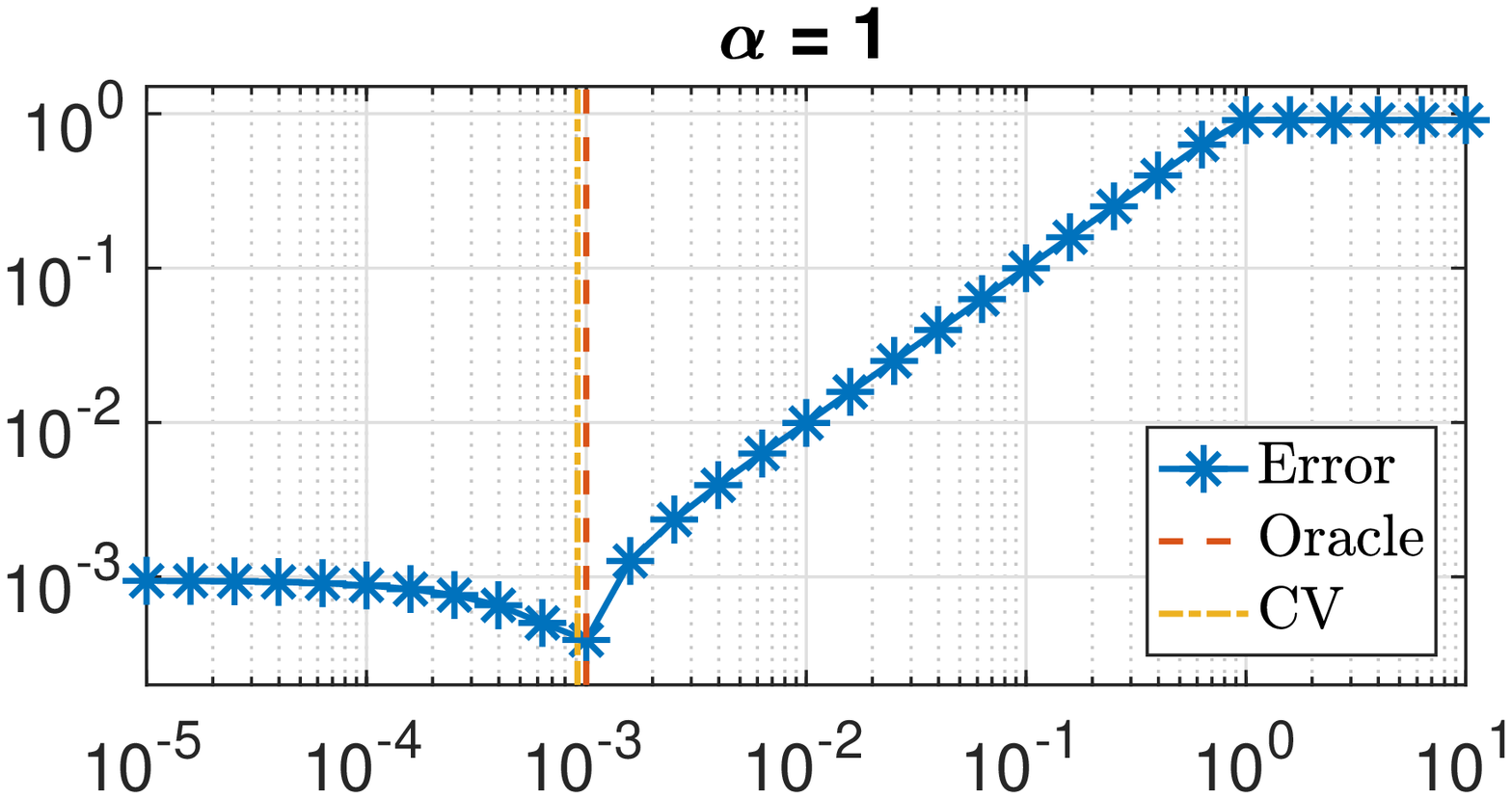}\\
\includegraphics[width=5.25cm]{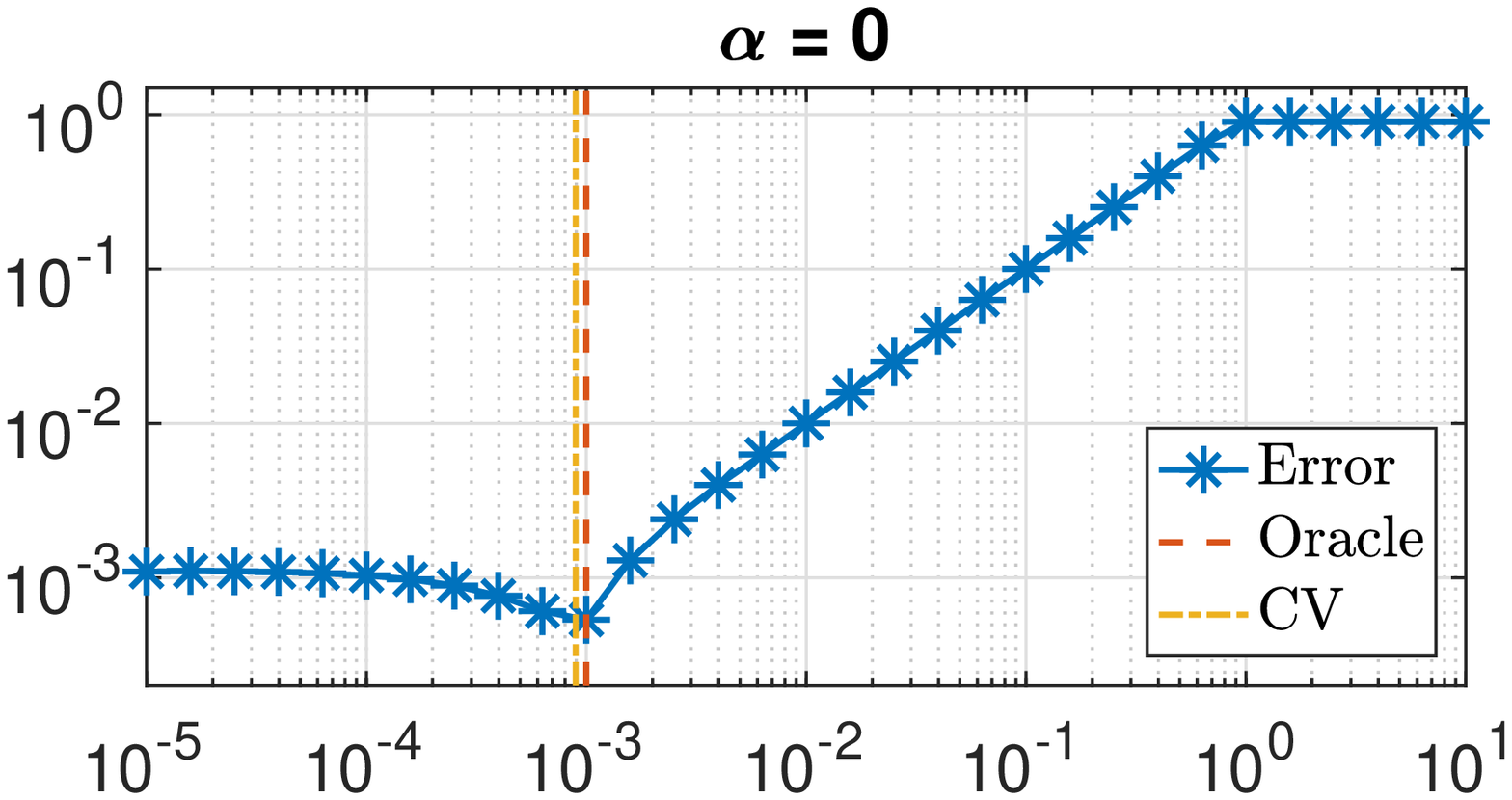} && 
\includegraphics[width=5.25cm]{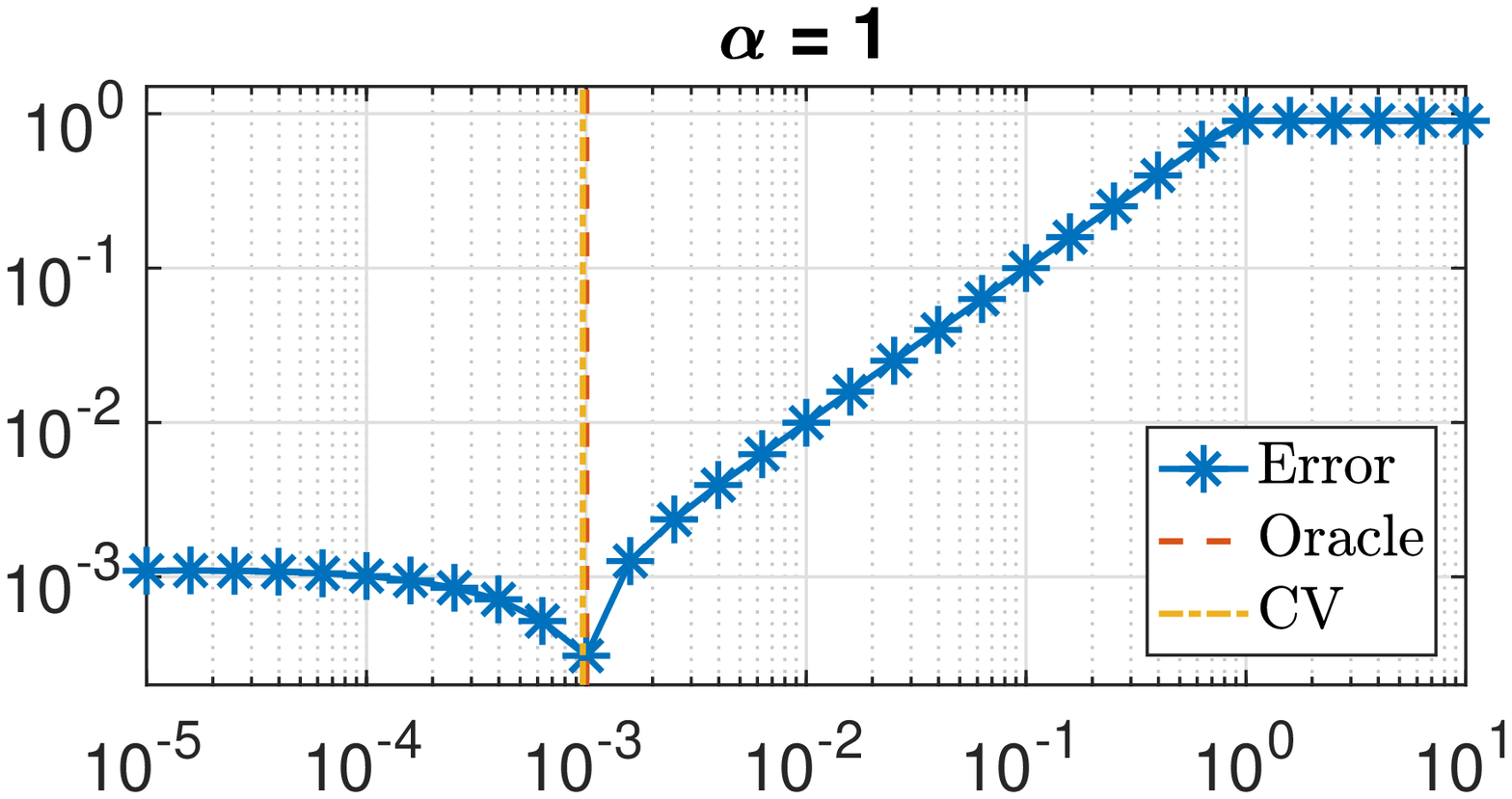}
\end{tabular}
\end{center}
\caption{\label{f:eta_vs_err_noise}Recovery error $\|f-\tilde{f}\|_{L^2_\nu}$ (averaged over 50 trials) against $\eta$, in the same setting as in Figure~\ref{f:m_vs_err_noise}. We use Chebyshev and Legendre polynomials in the top and bottom rows respectively. We consider $\eta = 10^{k}$, with $k$ belonging to a uniform grid of 31 points on the interval $[-5,1]$. The vertical lines represent the estimated values of $\eta$ (averaged over 50 trials) based on oracle least-squares (red dashed line) and cross validation (yellow dashed-dotted line). The weights are chosen as $\bm{w}=(\bm{u}_{\bm{i}}^\alpha)_{\bm{i}\in\Lambda}$, with $\alpha =0$ (left) and $\alpha =1$ (right).}
\end{figure}

\begin{table}
\centering
\begin{tabular}{|l|c|c|c|c|}
\hline
\multirow{2}{*}{$\alpha$} & \multicolumn{2}{|c|}{Chebyshev}      & \multicolumn{2}{|c|}{Legendre}       \\
\cline{2-5}
 & oracle & cross validation & oracle & cross validation \\ 
 \hline
0 & 1.0e-03 $\pm$ 7.2e-09 & 9.3e-04 $\pm$ 3.8e-04 & 1.0e-03 $\pm$ 3.7e-09 & 9.0e-04 $\pm$ 4.0e-04\\
1 & 1.0e-03 $\pm$ 4.9e-09 & 9.1e-04 $\pm$ 4.0e-04 & 1.0e-03 $\pm$ 3.6e-09 & 9.7e-04 $\pm$ 3.6e-04\\  
\hline
\end{tabular}
\caption{\label{t:mean_std_eta}Mean $\pm$ standard deviation for the values of $\eta$ estimated via oracle least-squares and cross validation over 50 trials in Figure~\ref{f:eta_vs_err_noise}.}
\end{table}

\section{Conclusions and challenges}\label{s:conclusions}

The concern of this chapter has been the emerging topic of compressed sensing for high-dimensional approximation.  As shown, smooth, multivariate functions are compressible in orthogonal polynomial bases.  Moreover, their coefficients have a certain form of structured sparsity corresponding to so-called lower sets.  The main result of this work is that such structure can be exploited via weighted $\ell^1$-norm regularizers.  Doing so leads to sample complexity estimates that are at most logarithmically dependent on the dimension $d$, thus mitigating the curse of dimensionality to a substantial extent.

As discussed in \S \ref{ss:literature}, this topic has garnered much interest over the last half a dozen years.  Yet challenges remain.  We conclude by highlighting a number of open problems in this area:

\textbf{Unbounded domains.}  We have considered only bounded hypercubes in this chapter.  The case of unbounded domains presents additional issues.  While Hermite polynomials (orthogonal on $\bbR$) have been considered in the case of unweighted $\ell^1$ minimization in \cite{GuoEtAlRandomizedQuad,HamptonDoostanCSPCE,NarayanJakemanZhouChristoffelLS}, the corresponding measurement conditions exhibit exponentially-large factors in either the dimension $d$ or degree $k$ of the (total degree) index space used.  It is unclear how to obtain dimension-independent measurement conditions in this setting, even for structured sparsity in lower sets.

\textbf{Sampling strategies.} Throughout we have considered sampling i.i.d.\ according to the orthogonality measure of the basis functions.  This is by no means the only choice, and various other sampling strategies have been considered in other works \cite{GuoEtAlRandomizedQuad,HamptonDoostanCSPCE,JakemanEtAlChristoffel,NarayanJakemanZhouChristoffelLS,NarayanZhouCCP,TangIaccarino,XuZhouSparseDeterministic}.  Empirically, several of these approaches are known to give some benefits.  However, it is not known how to design sampling strategies which lead to better measurement conditions than those given in Theorem \ref{t:unifquasiopt}.  A singular challenge is to design a sampling strategy for which $m$ need only scale linearly with $k$.  We note in passing that this has been achieved for the oracle least-squares estimator (recall \S \ref{ss:comparison}) \cite{MiglioratiCohenOptimal}.  However, it is not clear how to extend this approach to a compressed sensing framework.

\textbf{Alternatives to weighted $\ell^1$ minimization.}  As discussed in Remark \ref{r:structuredsparsity}, lower set structure is a type of structured sparsity model. We have used weighted $\ell^1$ minimization to promote such structure.  Yet other approaches may convey benefits.  Different, but related, types of structured sparsity have been exploited in the past using greedy or iterative algorithms \cite{BaranuikModelCS,BlumensathUnionSubspace,EldarDuarteCSReview,DuarteEldarStructuredCS}, or by designing appropriate convex regularizers \cite{TraonmilinGribonvalRIP}.  This remains an interesting problem for future work.

\textbf{Recovering Hilbert-valued functions.}
We have focused on compressed sensing-based polynomial approximation of high-dimensional functions whose coefficients belong to the complex domain $\mathbb{C}$.  However, an important problem in computational science, especially in the context of UQ and optimal control, involves the approximation of parametric PDEs.  Current 
compressed sensing techniques proposed in literature \cite{BouchotEtAlMultilevel,ChkifaDownwardsCS,DoostanOwhadiSparse,MathelinGallivanCSPDErandom,RS14,PengHamptonDoostantweighted,KarniadakisUQCS}  only approximate functionals of parameterized solutions, e.g.\ evaluation at a single spatial location, whereas a more robust approach should 
consider an $\ell_1$-regularized problem involving Hilbert-valued signals, i.e.\ signals where each coordinate is a function in a Hilbert space, which can provide a direct, global reconstruction of the solutions in the entire physical domain.  However, to achieve this goal new iterative minimization procedures as well as several theoretical concepts will need to be extended to the Hilbert space setting. 
The advantages of this approach over pointwise recovery with standard techniques will include: (i) for many parametric and stochastic model problems, global estimate of solutions in the physical domain is a quantity of interest; (ii) the recovery guarantees of this strategy can be derived from the decay of the polynomial coefficients in the relevant function space, which is well known in the existing theory; and (iii) the global reconstruction only assumes \textit{a priori} bounds of the tail expansion in energy norms, which are much more realistic than pointwise bounds.

\begin{acknowledgement}
The first and second authors acknowledge the support of the Alfred P. Sloan Foundation and the Natural Sciences and Engineering Research Council of Canada through grant 611675.
The second author acknowledges the Postdoctoral Training Center in Stochastics of the Pacific Institute for the Mathematical Sciences for the support.
The third author acknowledges support by: the U.S. Defense Advanced Research Projects Agency, Defense Sciences Office under contract and award numbers HR0011619523 and 1868-A017-15; the U.S. Department of Energy, Office of Science, Office of Advanced Scientific Computing Research, Applied Mathematics program under contract number ERKJ259 and ERKJ314, and; the Laboratory Directed Research and Development program at the Oak Ridge National Laboratory, which is operated by UT-Battelle, LLC., for the U.S. Department of Energy under Contract DE-AC05-00OR22725.
\end{acknowledgement}

\bibliographystyle{plain}
\bibliography{MatheonRefs}

\end{document}